\documentclass[hidelinks,onefignum,onetabnum]{siamart251216}

\usepackage[sort&compress,numbers]{natbib}
\usepackage{amsfonts}
\usepackage{amssymb}
\usepackage{booktabs}

\ifpdf
  \DeclareGraphicsExtensions{.eps,.pdf,.png,.jpg}
\else
  \DeclareGraphicsExtensions{.eps}
\fi

\newsiamremark{remark}{Remark}
\newsiamremark{example}{Example}
\newsiamremark{assumption}{Assumption}
\newsiamthm{claim}{Claim}

\newcommand{\one}{\mathbf{1}}

\newcommand{\Rnum}{\mathbb{R}}

\newcommand{\KL}{\mathrm{KL}}

\newcommand{\cM}{\mathcal{M}}
\newcommand{\cP}{\mathcal{P}}
\newcommand{\cX}{\mathcal{X}}

\newcommand{\WFR}{\mathrm{WFR}}

\newcommand{\dd}{\,\mathrm{d}}

\headers{Discrete WFR}
        {J. Ying, X. Wang, P. Zhou, and L. Zhang}

\title{Wasserstein--Fisher--Rao Reconstruction of Discrete Space Population Dynamics%
\thanks{Submitted to the editors [DATE TO BE ADDED].
\funding{This work was supported by the National Natural Science Foundation of China (NSFC No. 12225102 to L.Z.,
8206100646 to P.Z., T2321001 to P.Z. \& L.Z., and 12288101
to P.Z. \& L.Z.) and the National Key R\&D Program of
China No.2024YFA0919500 to L.Z.}}}

\author{Junda Ying\hspace{0.4em}\thanks{Beijing International Center for Mathematical
Research, Peking University (\email{yingjunda@stu.pku.edu.cn}).}
\and Xinyu Wang\hspace{0.4em}\thanks{LMAM and School of Mathematical Sciences, Peking University. (\email{xinyuwang@pku.edu.cn}).}
\and Peijie Zhou\thanks{Corresponding author.
Center for Quantitative Biology and Center for Machine
Learning Research, Peking University; National Engineering Laboratory for Big
Data Analysis and Applications, Beijing; AI for Science Institute, Beijing
(\email{pjzhou@pku.edu.cn}).}
\and Lei Zhang\thanks{Corresponding author.
Beijing International Center for Mathematical Research,
School of Mathematical Sciences, Center for Quantitative Biology, Center for
Machine Learning Research, and Institute for Artificial Intelligence, Peking
University (\email{zhangl@math.pku.edu.cn}).}}

\ifpdf
\hypersetup{
  pdftitle={Wasserstein--Fisher--Rao Reconstruction of Discrete space Population Dynamics},
  pdfauthor={Xinyu Wang, Junda Ying, Peijie Zhou, and Lei Zhang}
}
\fi

\begin{document}

\maketitle
\begin{abstract}
Reconstructing the continuous time dynamics of cell state evolution from sparse snapshot observations can be formulated as an unbalanced optimal transport problem. To study such dynamics at a coarse grained cell type level, unbalanced optimal transport must be extended to discrete state spaces. Inspired by discrete flow matching, we extended the dynamic Wasserstein–Fisher–Rao (WFR) problem to discrete spaces, and derived a closed form for its solution. For particularly large state spaces or more general models, we proposed an unbalanced discrete flow matching as a numerical approximation method, which provides a Markovian reconstruction of the optimal WFR flow. Our discrete WFR framework also stays stable across different coarse grain resolution. Finally, numerical experiments demonstrated the effectiveness of the proposed approach, and applications to real biological data supported the plausibility of the reconstructed dynamics.
\end{abstract}



\section{Introduction}\label{intro}
Mathematical models of cell population dynamics play an important role in understanding embryonic development and cellular differentiation. Structured population models and multicompartment models have been developed to describe these processes by accounting for cell state transitions, proliferation, and death \cite{doumic2011structured,zhang2022continuous}. With the increasing availability of single-cell transcriptomic (scRNA-seq) data, recent studies have sought to reconstruct cellular dynamics directly from single-cell level observations \cite{jiang2022dynamic,jiang2024physics,lavenant2024toward,SDEvelo}. However, due to the destructive nature of sequencing, such measurements typically provide unpaired snapshots rather than longitudinal trajectories of individual cells. A natural mathematical question is therefore how to reconstruct continuous time unbalanced dynamics that describe both cell state transitions and population growth from a finite collection of sparsely sampled, unpaired snapshots.

Mathematically, if each cell is represented as a vector in \(\mathbb{R}^d\), the snapshots observed at different time points can be regarded as a collection of measures \(\mu_1,\ldots,\mu_K\) with possibly different total masses. Dynamical reconstruction then seeks a time dependent velocity field \(v_t(x)\) and a time dependent growth rate \(g_t(x)\) on \(\mathbb{R}^d\) such that the associated unbalanced continuity equation
\[
\partial_t\mu_t+\nabla\cdot(\mu_t v_t)=g_t\mu_t
\]
interpolates these observed measures. Since infinitely many dynamical systems may satisfy these interpolation constraints, optimal transport (OT) provides a natural variational principle for selecting interpolating dynamics. Its dynamic formulation characterizes mass transport through a least-action principle \cite{benamou2000computational}, while unbalanced extensions such as the Wasserstein--Fisher--Rao (WFR) metric further incorporate mass creation and destruction \cite{chizat2015interpolating,liero2016optimal}. These mathematical frameworks have motivated OT-based approaches to cellular trajectory inference, growth modeling, and stochastic dynamical reconstruction \cite{waddington-ot,TIGON,DeepRUOT,cang2025synchronized,cang2026synchronization,liu2026learning} with nice physical interpretation that physical and biological processes tend to save energy. More recently, flow matching has provided conditional path based approaches to learning balanced and unbalanced transport dynamics \cite{cfm_tong,wang2025joint,wfrfm2025}.

As scRNA-seq datasets increasingly contain millions of cells \cite{qiu2024single}, reconstructing dynamics at single-cell resolution poses substantial computational challenges, despite recent advances in numerical methods for high-dimensional dynamic OT \cite{wan2023scalable}. Moreover, many biological questions concern population abundances, transition rates, and growth rates associated with a finite collection of cell types rather than individual trajectories in the high-dimensional gene expression space. This motivates the study of cellular dynamics at a coarse-grained cell type level \cite{jiang2022dynamic,doumic2011structured,zhang2022continuous}. At this level of description, the state space changes from \(\mathbb{R}^d\) to a finite set \(\mathcal{X}\), and the unbalanced continuity equation is replaced by an unbalanced Kolmogorov forward equation:
\begin{equation}
\label{eq:Markov dynamics}
\partial_t p_t(x)
=
\sum_{y\neq x}p_t(y)q_t(y,x)
-
p_t(x)\sum_{y\neq x}q_t(x,y)
+
p_t(x)g_t(x),
\qquad x\in\mathcal X.
\end{equation}
The corresponding dynamics inference problem is then to find transition rates \(q_t\) and growth rate \(g_t\) that interpolate the observed state-level measures while minimizing an appropriate energy. Inspired by optimal transport and Wasserstein–Fisher–Rao (WFR) formulations in continuous spaces \cite{benamou2000computational,chizat2015interpolating}, these considerations motivate the study of dynamic WFR directly on discrete spaces.

Dynamic OT on discrete spaces has been extensively studied through the extension of Wasserstein geometry and gradient flow structures to graphs and finite Markov chains \cite{ChowLiZhou2018Entropy,ChowHuangLiZhou2012FokkerPlanck, Maas2011Gradient,Zhou2021OptimalTransport,gladbach2020scaling}. Related studies have also developed discretizations of unbalanced OT and their convergence theories \cite{xiong2024convergence,li2024convergence}, as well as transport metrics for nonnegative measures on finite Markov chains \cite{mao2026wasserstein}. Rather than defining a transport action through a prescribed graph-dependent mobility, we construct a dynamic realization of discrete space WFR through endpoint-conditioned transport and growth paths, inspired by discrete flow matching \cite{gat2024discrete}. The resulting measure flow and its weighted Markovian reconstruction \eqref{eq:Markov dynamics} can be explicitly obtained from an optimal static semi-coupling without further numerical approximation. Given proper conditions, our dynamic discrete WFR measure flow also stays consistent under state aggregation, therefore robust to cell type annotations at different resolutions.

The contribution is the connection between three structures: discrete static WFR semi-couplings, dynamic WFR measure flows, and its Markovian reconstruction. First, we develop a dynamic realization of discrete WFR based on travelling Dirac paths, and prove the equivalence between it and the static discrete WFR, as well as an optimal entropy transport version which can be efficiently solved. Second, we derive an explicit weighted Markovian representation of the resulting measure flow, with unbalanced discrete flow matching as an extension when direct marginalization is impractical. Third, we establish consistency and stability for discrete WFR action and measure flow under state aggregation. 

The paper is organized as follows. Section \ref{sec: WFR theory} formulates the dynamic WFR problem on discrete spaces and investigates its basic mathematical properties, including travelling Dirac solutions, static and dynamic form, and the metric property. Section \ref{sec:unbalanced_dfm} constructs the dynamic WFR measure flow, its explicit Markovian reconstruction, and a more general unbalanced discrete flow matching framework. Section \ref{sec:coarse_consistency} studies consistency under state aggregation. Section~\ref{sec:experiments} reports numerical validations and biological applications. Section~\ref{conclusion} summarizes the work.

\section{Discrete Wasserstein-Fisher-Rao distance}
\label{sec: WFR theory}

We develop a dynamic Wasserstein--Fisher--Rao model on a discrete space. The construction follows the same logic as the continuous theory. We first recall how balanced optimal transport passes from the static Kantorovich problem to the dynamic Benamou--Brenier formulation, and how WFR augments the dynamic equation with growth before returning to a static semi-coupling problem. We then repeat this chain on a finite metric space. We first define the static and dynamic OT on discrete space, and then add growth penalty to define a dynamic WFR problem on discrete space. Finally, we derive the static discrete WFR problem, and study its several properties and an equivalent optimal entropy transport form which can be efficiently solved.

\subsection{Continuous OT and WFR}\label{subsec:continuous-guide}

Let $\Omega\subset\Rnum^d$ and let $\mu_0,\mu_1\in\cM_+(\Omega)$ be densities of absolutely continuous measures with unit total mass. The static Kantorovich formulation of quadratic optimal transport is
\begin{equation}\label{eq:cont-static-ot}
    \mathcal{W}_2^2(\mu_0,\mu_1)
    =
    \inf_{\pi\in\Pi(\mu_0,\mu_1)}
    \int_{\Omega\times\Omega}|x-y|^2\,\pi(x,y)\dd x\dd y,
\end{equation}
where $\Pi(\mu_0,\mu_1)$ is the set of couplings with marginals $\mu_0$ and $\mu_1$ \cite{kantorovich1942transfer,villani2009optimal}. The dynamic Benamou--Brenier formulation writes the same quantity as
\begin{equation}\label{eq:cont-bb}
    \mathcal{W}_2^2(\mu_0,\mu_1)
    =
    \inf_{\mu_t,v_t}
    \int_0^1\int_{\Omega}|v_t(x)|^2\mu_t(x)\dd x\dd t,
    \qquad
    \partial_t\mu_t+\nabla\cdot(\mu_t v_t)=0,
\end{equation}
with endpoint conditions $\mu_{t=0}=\mu_0$ and $\mu_{t=1}=\mu_1$ \cite{benamou2000computational}. This dynamic viewpoint identifies optimal transport as a least-action problem with a conservation law constraint.

Wasserstein--Fisher--Rao (WFR) \cite{chizat2015interpolating,chizat2018unbalanced,liero2016optimal,liero2018optimal} removes the unit mass constraint by adding a growth penalty. In the normalization used in this paper, the continuous dynamic WFR problem has the formal action
\begin{equation}\label{eq:cont-wfr-dynamic}
    \WFR_\delta^2(\mu_0,\mu_1)
    =
    \inf_{\mu_t,v_t,g_t}
    \int_0^1\int_{\Omega}
    \left(|v_t(x)|^2+\delta^2 g_t(x)^2\right)\mu_t(x)\dd x\dd t,
\end{equation}
subject to
\begin{equation*}
    \partial_t\mu_t+\nabla\cdot(\mu_t v_t)=g_t\mu_t,
    \qquad
    \mu_{t=0}=\mu_0,\quad \mu_{t=1}=\mu_1.
\end{equation*}
The velocity field $v_t$ transports existing mass, while $g_t$ creates or annihilates mass. By direct calculation, the WFR cost between two Dirac measures \(a\delta_x,b\delta_y\) is
\begin{equation}\label{eq:cont-dirac-cost}
    4\delta^2
    \left[
        a+b-2\sqrt{ab}\,
        \overline{\cos}\left(\frac{|x-y|}{2\delta}\right)
    \right],
\end{equation}
where $a$ and $b$ are the masses attached to the endpoint pair $(x,y)$, \(\overline{\cos}(x)=\cos(x\land\frac{\pi}{2})\). The corresponding optimal trajectory of the Dirac measure is called the travelling Dirac. The key structural point is that the travelling Dirac cost is convex and positively one-homogeneous when \(|x-y|<\pi\delta\), hence sublinear in the two endpoint masses \((a,b)\). This sublinearity allows multiple travelling Dirac components with the same endpoints to be merged without increasing cost, which reduces the dynamic WFR problem (\eqref{eq:cont-wfr-dynamic}) into a static form
\begin{equation}
\begin{aligned}
\WFR_\delta^2(\mu_0,\mu_1)
&=4\delta^2\inf_{(\gamma_0,\gamma_1)}\int_{\Omega^2}
\Biggl(\gamma_0(x,y)+\gamma_1(x,y)\\
&\qquad
-2\sqrt{\gamma_0(x,y)\gamma_1(x,y)}
\overline{\cos}\left(\frac{|x-y|}{2\delta}\right)\Biggr)
\mathrm{d}x\,\mathrm{d}y
\end{aligned}
\end{equation}
subject to 
\begin{equation*}
\begin{aligned}
(\gamma_0,\gamma_1)\in\Gamma(\mu_0,\mu_1)
=\Biggl\{(\gamma_0,\gamma_1)&\in
\left(\mathcal{M}_+(\Omega^2)\right)^2:\\
&\int_\Omega\gamma_0(x,y)\mathrm{d}y=\mu_0(x),\quad
\int_\Omega\gamma_1(x,y)\mathrm{d}x=\mu_1(y)\Biggr\}.
\end{aligned}
\end{equation*}
where \((\gamma_0,\gamma_1)\) is called the semi-coupling as it is the unbalanced version of the optimal transport coupling. This formulation is said to be static since there is no time variable \(t\) in it. It can be further reduced to an optimal entropy transport problem
\begin{equation}
\label{OET}
\begin{aligned}
\WFR_{\delta}^2(\mu_0,\mu_1) &= 4\delta^2\inf_{\gamma\in\mathcal{M}_+(\Omega^2)} \{\int_{\Omega^2}  \, -2\operatorname{ln}\operatorname{\overline{\cos}}\left(\frac{|x-y|}{2\delta}\right)\gamma(x,y)\mathrm{d}x\mathrm{d}y\\&+\operatorname{KL}(\int_\Omega\gamma(x,y)\mathrm{d}y\Vert \mu_0(x))+\operatorname{KL}(\int_\Omega\gamma(x,y)\mathrm{d}x\Vert \mu_1(y))\}
\end{aligned}
\end{equation}
which can be solved efficiently. Recent research has also utilized flow matching for efficiently solving the dynamic WFR problem \cite{wfrfm2025}. The WFR framework defines a metric on \(\mathcal{M}_+(\Omega)\), and its dynamic form models the unbalanced dynamics in \(\mathbb{R}^d\). In this paper, we extend the framework to graphs to model unbalanced dynamics on discrete spaces.

\subsection{Dynamic discrete OT}\label{subsec:discrete-ot}
Let \(\mathcal{G}=(\mathcal{X},C)\) be a finite graph where \(\mathcal{X}=\{x_1,\ldots,x_N\}\) is a finite set, and $C_{ij} \in [0,+\infty]$ are weighted edges. We write \(\cM_{+}(\cX)=\left\{\rho=(\rho_1,\ldots,\rho_N)\in \Rnum_{\ge0}^N\right\}\) for the set of nonnegative measures over $\cX$. The total mass is denoted by \(|\rho|=\sum_{i=1}^N \rho_i\). \(\mathcal{P}(\mathcal{X})\) stands for the set of probability measures over \(\mathcal{X}\). 

For $\rho,\eta\in\cP(\cX)$, the static OT cost between them is defined as 
\begin{equation}\label{eq:discrete-static-ot}
    \mathrm{OT}(\rho,\eta)
    =
    \min_{\pi\in\cP(\cX^2)}
    \left\{
        \sum_{i,j=1}^N C_{ij}\pi_{ij}:
        \sum_j\pi_{ij}=\rho_i,\quad
        \sum_i\pi_{ij}=\eta_j
    \right\}.
\end{equation}
As a well known result, when \(C_{ij}\) is a metric on \(\mathcal{X}\), OT is also a metric on \(\cP(\cX)\). For more general cases, we have the following proposition.

\begin{proposition}[Metric property of OT on graph]\label{prop:ot metric}
For \(p\ge1\), if \(C_{ij}^{\frac{1}{p}}\) defines a metric on \(\cX\), \(\cP(\cX)\) is a metric space equipped with the metric \(\mathrm{OT}^{\frac{1}{p}}\). When \(p=1\), \(\cP(\cX)\) is further a geodesic space under this metric, with geodesic connecting \(\alpha,\beta\) defined as \(\gamma(t)=(1-t)\alpha+t\beta\).
\end{proposition}
\begin{proof}
The proof is provided in Supplementary Material~\ref{app:proof-ot-metric}.    
\end{proof}

Inspired by discrete flow matching \cite{campbell2024generative,gat2024discrete}, the dynamic form is defined by assigning each pair of mass transport \(\pi_{ij}\) a continuous time probability flow. For a coupling $\pi$ and progress variables $k_{ij}:[0,1]\to[0,1]$ with $k_{ij}(0)=0$ and $k_{ij}(1)=1$, define
\begin{equation}\label{eq:discrete-ot-path}
    \pi_{ij}
    \left[
        (1-k_{ij}(t))\delta_{x_i}
        +
        k_{ij}(t)\delta_{x_j}
    \right].
\end{equation}
as the probability flow transporting mass \(\pi_{ij}\) from \(x_i\) to \(x_j\). Since \(C_{ij}\) is at the position of the squared Euclidean distance in the continuous case, the product of \(C_{ij}^{\frac{1}{2}}\) and the derivative of the progress variable \(\dot{k}_{ij}\) measures the speed of transport. Hence, the corresponding dynamic action is defined as an analog of the kinetic action in Benamou--Brenier form 
\begin{equation}\label{eq:discrete-dynamic-ot}
    \mathcal A_{\mathrm{OT}}[\pi,k]
    =
    \int_0^1
    \sum_{i,j}\pi_{ij}C_{ij}\dot k_{ij}(t)^2\,\dd t .
\end{equation}
subject to
\[
\sum_j\pi_{ij}=\rho_i, \quad\sum_i\pi_{ij}=\eta_j,\quad k_{ij}(0)=0,\quad k_{ij}(1)=1
\]

\begin{theorem}[Static-dynamic equivalence for OT]\label{prop:finite-ot-equivalence}
For $\rho,\eta\in\cP(\cX)$,
\begin{equation}\label{eq:finite-ot-equivalence}
    \mathrm{OT}(\rho,\eta)
    =
    \inf_{\pi,k}\mathcal A_{\mathrm{OT}}[\pi,k],
\end{equation}
\end{theorem}

\begin{proof}
For fixed $\pi$, the action separates over pairs. For each active pair,
\begin{equation*}
    \int_0^1 \dot k_{ij}(t)^2\,\dd t
    \geq
    \left(\int_0^1\dot k_{ij}(t)\,\dd t\right)^2
    =
    1
\end{equation*}
by Cauchy's inequality, with equality for $k_{ij}(t)=t$. Hence, the least dynamic action for fixed $\pi$ is $\sum_{i,j}C_{ij}\pi_{ij}$. Minimizing over $\pi$ gives \eqref{eq:discrete-static-ot}.
\end{proof}

\begin{remark} Although inspired by viewing \(C_{ij}\) as an analog of the squared Euclidean distance, the definition of the dynamic action does not require \(C_{ij}^{\frac{1}{2}}\) to be a metric. 
\end{remark}
When \(C_{ij}\) is a metric, the induced continuous probability flow
\[
\sum_j\left(\pi_{ij}(1-t)+\pi_{ji}t\right)
=(1-t)\rho_i+t\eta_i
\]
recovers the geodesic in \(\cP(\cX)\). These results indicate that our defined dynamics formulation is reasonable.
\subsection{Dynamic discrete WFR}\label{subsec:dynamic WFR}

In this section, we extend the above to unbalanced scenarios. We first consider a pair of weighted Dirac masses $m_0\delta_{x_i}$ and $m_1\delta_{x_j}$ with $m_0,m_1\ge0$. A path is represented by a scalar mass $m(t)\ge0$ and a progress variable $k(t)\in[0,1]$. Formally, $k(t)=0$ corresponds to the source point $x_i$ and $k(t)=1$ corresponds to the target point $x_j$. Adding a quadratic growth penalty scaled by a parameter \(\delta\) to the OT action, the unbalanced action is defined as
\begin{equation}\label{eq:dirac-action}
\begin{aligned}
    \mathcal{A}_{ij}[m,k]
    &=\int_0^1
    m_{ij}(t)\left[
        C_{ij}\dot{k}_{ij}(t)^2
        +\delta^2\frac{\dot{m}_{ij}(t)^2}{m_{ij}(t)^2}
    \right]\dd t\\
    &=\int_0^1
    \left[
        m_{ij}(t)C_{ij}\dot{k}_{ij}(t)^2
        +\delta^2\frac{\dot{m}_{ij}(t)^2}{m_{ij}(t)}
    \right]\dd t
\end{aligned}
\end{equation}
with the convention that the growth term is $+\infty$ unless $m$ is absolutely continuous and $\dot m=0$ a.e. on $\{m=0\}$. The endpoints are
\begin{equation*}
    m_{ij}(0)=m_0,\qquad m_{ij}(1)=m_1,\qquad k_{ij}(0)=0,\qquad k_{ij}(1)=1.
\end{equation*}

\begin{definition}[Dirac-pair WFR cost on graph]\label{def:diracwfr}
For $m_0,m_1\geq0$ and $i,j\in\{1,\ldots,N\}$, define
\begin{equation}\label{eq:dirac-def}
    \WFR_{\delta}^2(m_0\delta_{x_i},m_1\delta_{x_j})
    =
    \inf_{m,k}\mathcal{A}_{ij}[m,k],
\end{equation}
where the infimum is taken over all admissible paths satisfying the above endpoint conditions.
\end{definition}

\begin{proposition}[Travelling Dirac on graph]\label{prop:dirac}
For any $m_0,m_1>0$, and \(C_{ij}<\pi^2\delta^2\),
\begin{equation}\label{eq:dirac-closed}
    \WFR_{\delta}^2(m_0\delta_{x_i},m_1\delta_{x_j})
    =
    4\delta^2
    \left[
        m_0+m_1
        -2\sqrt{m_0m_1}\,
        \cos\left(\frac{\sqrt{C_{ij}}}{2\delta}\right)
    \right].
\end{equation}
The corresponding travelling Dirac path is
\[m_{ij}(t)=At^2+Bt+m_0,\quad k_{ij}(t)=\omega\int_0^t\frac{\mathrm{d}s}{m_{ij}(s)}\]
where 
\begin{equation}
\label{eq:travelling Dirac}
\begin{aligned}
&A=m_0+m_1-\sqrt{4m_0m_1-\frac{\omega^2C_{ij}}{\delta^2}}, \quad B=-2m_0+\sqrt{4m_0m_1-\frac{\omega^2C_{ij}}{\delta^2}}\\
&\omega=2\delta\tau\sqrt{\frac{m_0m_1}{C_{ij}(1+\tau^2)}},\quad \tau=\tan\frac{\sqrt{C_{ij}}}{2\delta}  
\end{aligned}
\end{equation}
\end{proposition}
\begin{proof}

Denote \(\theta_{ij}:=\frac{\sqrt{C_{ij}}}{2\delta}.\) When \(C_{ij}<\pi^2\delta^2\), we have \(0\le \theta_{ij}<\pi/2\). Consider the travelling paths with both transportation and birth-death
\[
    m(t)[(1-k(t))\delta_{x_i}+k(t)\delta_{x_j}],
\]
where \(k(0)=0\), \(k(1)=1\), \(m(0)=m_0\), and \(m(1)=m_1\). Let \(r(t)=\sqrt{m(t)},\theta(t)=\frac{\sqrt{C_{ij}}}{2\delta}k(t).\) Then the travelling Dirac action becomes
\begin{equation}
\begin{aligned}
    \int_0^1\left[C_{ij}m(t)\dot k(t)^2
    +\delta^2\frac{\dot m(t)^2}{m(t)}\right]\mathrm dt
    =
    4\delta^2\int_0^1\left[\dot r(t)^2+r(t)^2\dot\theta(t)^2\right]\mathrm dt .
\end{aligned}
\end{equation}
Introduce the polar coordinates \(z(t)=r(t)e^{\mathrm i\theta(t)}.\) Then \(|\dot z(t)|^2=\dot r(t)^2+r(t)^2\dot\theta(t)^2.\) Therefore, the minimization among travelling Dirac paths is reduced to the Euclidean shortest path problem \(\min_z 4\delta^2\int_0^1|\dot z(t)|^2\mathrm dt\) with boundary conditions \(z(0)=\sqrt{m_0},z(1)=\sqrt{m_1}e^{\mathrm i\theta_{ij}}.\) The minimizer is the straight line
\[
    z^*(t)=(1-t)\sqrt{m_0}+t\sqrt{m_1}e^{\mathrm i\theta_{ij}}.
\]
Thus
\begin{equation}
\begin{aligned}
    \min 4\delta^2\int_0^1|\dot z(t)|^2\mathrm dt
    =4\delta^2\left|\sqrt{m_1}e^{\mathrm i\theta_{ij}}-\sqrt{m_0}\right|^2
    =4\delta^2\left[m_0+m_1-2\sqrt{m_0m_1}\cos\theta_{ij}\right].
\end{aligned}
\end{equation}

Finally, we derive the explicit form of the path. From the straight line formula,
\begin{equation}
\label{eq:travelling Dirac mass}
m_{ij}(t)=|z^*(t)|^2=(1-t)^2m_0+t^2m_1+2t(1-t)\sqrt{m_0m_1}\cos\theta_{ij}.   
\end{equation}
    
Thus, \(m_{ij}(t)=At^2+Bt+m_0\) where
\[
    A=m_0+m_1-2\sqrt{m_0m_1}\cos\theta_{ij},\qquad
    B=-2m_0+2\sqrt{m_0m_1}\cos\theta_{ij}.
\]
Let \(\tau=\tan\theta_{ij},\omega=2\delta\tau\sqrt{\frac{m_0m_1}{C_{ij}(1+\tau^2)}}.\) Since \(0\le\theta_{ij}<\pi/2\),\(\omega=\frac{2\delta\sqrt{m_0m_1}\sin\theta_{ij}}{\sqrt{C_{ij}}}.\) Moreover, \(\sqrt{4m_0m_1-\frac{\omega^2C_{ij}}{\delta^2}}=2\sqrt{m_0m_1}\cos\theta_{ij}.\) Hence \(A\) and \(B\) can be equivalently written as
\[
    A=m_0+m_1-\sqrt{4m_0m_1-\frac{\omega^2C_{ij}}{\delta^2}},
    \qquad
    B=-2m_0+\sqrt{4m_0m_1-\frac{\omega^2C_{ij}}{\delta^2}}.
\]
The angular velocity of the straight line satisfies
\[
\dot k_{ij}(t)=\frac{\omega}{m_{ij}(t)}.
\]
Consequently,
\[
k_{ij}(t)=\omega\int_0^t\frac{\mathrm ds}{m_{ij}(s)}.
\]

\end{proof}

\begin{remark}[Zero mass endpoint limit]
For paths with zero mass endpoints, i.e., $m_0=0$ or $m_1=0$, the optimal path is the limit of the paths with positive endpoints.
$$
\left\{
\begin{aligned}
&m_{ij}(t)=m_1t^2,\quad k_{ij}(t)=\mathrm{1}_{t>0}(t),\qquad m_0=0.\\
&m_{ij}(t)=m_0(1-t)^2,\quad k_{ij}(t)=\mathrm{1}_{t\ge1}(t),\qquad m_1=0.   
\end{aligned}
\right .
$$
\end{remark}

\begin{remark}[Infinite growth penalty limit]
If $m_0=m_1=m$ and $\delta\gg1$, then
\begin{equation*}
\begin{aligned}
\WFR_{\delta}^2(m\delta_{x_i},m\delta_{x_j})
&=8\delta^2m\left[1-\cos\left(\frac{\sqrt{C_{ij}}}{2\delta}\right)\right]\\
&=mC_{ij}+O(\delta^{-2})\to mC_{ij}
\quad(\delta\to+\infty).
\end{aligned}
\end{equation*}
Thus, the WFR cost recovers the balanced OT cost as the growth penalty \(\delta\) tends to infinity.
\end{remark}

\begin{proposition}[Sublinearity of the travelling Dirac cost]\label{prop:dirac-sublinear}
Fix an endpoint pair $(i,j)$ and set
\begin{equation*}
    D_{ij}(a,b)
    =
    4\delta^2
    \left[
        a+b-2\sqrt{ab}\,
        \cos\left(\frac{\sqrt{C_{ij}}}{2\delta}\right)
    \right],
    \qquad a,b\geq0.
\end{equation*}
Then $D_{ij}$ is positively one-homogeneous and sublinear:
\begin{equation}\label{eq:dirac-subadditive}
    D_{ij}\left(\sum_{\ell}a_{\ell},\sum_{\ell}b_{\ell}\right)
    \leq
    \sum_{\ell}D_{ij}(a_{\ell},b_{\ell})
\end{equation}
for all finite nonnegative $(a_\ell,b_\ell)$.
\end{proposition}

\begin{proof}
Positive one-homogeneity is immediate from the formula. Write
\begin{equation*}
\zeta_{ij}=\cos\left(\frac{\sqrt{C_{ij}}}{2\delta}\right)\in[0,1],
    \qquad
    A=\sum_\ell a_\ell,\qquad B=\sum_\ell b_\ell .
\end{equation*}
By Cauchy's inequality,
\begin{equation*}
    \sum_\ell \sqrt{a_\ell b_\ell}
    \leq
    \sqrt{\left(\sum_\ell a_\ell\right)
    \left(\sum_\ell b_\ell\right)}
    =
    \sqrt{AB}.
\end{equation*}
Since $\zeta_{ij}\geq0$,
\begin{equation*}
\begin{split}
    \sum_\ell D_{ij}(a_\ell,b_\ell)
    &=
    4\delta^2
    \left[
        A+B-2\zeta_{ij}\sum_\ell\sqrt{a_\ell b_\ell}
    \right]\\
    &\geq
    4\delta^2
    \left[
        A+B-2\zeta_{ij}\sqrt{AB}
    \right]
    =
    D_{ij}(A,B).
\end{split}
\end{equation*}
This proves \eqref{eq:dirac-subadditive}.
\end{proof}

With sublinearity, multiple unbalanced transport paths from \(x_i\) to \(x_j\) tend to merge into one to minimize the total transport cost. Hence, it is reasonable to define the dynamic WFR action as following.
\begin{definition}[Dynamic WFR on graph]
For \(\rho,\eta\in\cM_+(\cX)\), define
\begin{equation}\label{eq:discrete-dynamic-wfr}
    \WFR_{\delta}^2(\rho,\eta)
    =
    \inf_{m,k}\int_0^1
    \sum_{i,j}\left(m_{ij}(t)C_{ij}\dot{k}_{ij}(t)^2
        +\delta^2\frac{\dot{m}_{ij}(t)^2}{m_{ij}(t)}\right)\dd t .
\end{equation}
subject to
\[
\sum_j m_{ij}(0)=\rho_i, \quad\sum_i m_{ij}(1)=\eta_j,\quad k_{ij}(0)=0,\quad k_{ij}(1)=1
\]  
\end{definition}

\subsection{Discrete static WFR}\label{subsec:static}
As a discrete version of the semi-coupling, let $a=(a_{ij})$ and $b=(b_{ij})$ be two nonnegative matrices satisfying
\begin{equation}\label{eq:ab-marginals}
    \sum_{j=1}^N a_{ij}=\rho_i,\qquad
    \sum_{i=1}^N b_{ij}=\eta_j.
\end{equation}
We denote this semi-coupling set by
\[
\Gamma(\rho,\eta)
:=\left\{(a,b)\geq0:
\sum_j a_{ij}=\rho_i,\quad
\sum_i b_{ij}=\eta_j\right\}.
\]
The quantity $a_{ij}$ is the amount of mass initially assigned to the pair $(i,j)$, while $b_{ij}$ is the amount of mass finally assigned to the same pair. In unbalanced settings, \(a_{ij},b_{ij}\) can be different. By (\ref{eq:dirac-closed}), when \(\max_{i,j}{C_{ij}}<\pi^2\delta^2\), and given \(m_{ij}(0)=a_{ij},m_{ij}(1)=b_{ij}\) fixed, the dynamic WFR (\ref{eq:discrete-dynamic-wfr}) can be minimized by first minimizing over \(k_{ij},m_{ij}\), then \(a_{ij},b_{ij}\), which yields a static form for WFR on graph.

\begin{definition}[Static WFR on graph]\label{def:staticwfr}
For $\rho,\eta\in\cM_+(\cX)$, \(\max_{i,j}{C_{ij}}<\pi^2\delta^2\), define
\begin{equation}\label{eq:static-wfr}
    \WFR_{\delta}^2(\rho,\eta)
    =
    4\delta^2
    \inf_{(a,b)\in\Gamma(\rho,\eta)}
    \sum_{i,j=1}^N
    \left[
        a_{ij}+b_{ij}
        -2\sqrt{a_{ij}b_{ij}}\,
        \cos\left(\frac{\sqrt{C_{ij}}}{2\delta}\right)
    \right].
\end{equation}
\end{definition}

Several basic properties were investigated as following.
\begin{proposition}[Existence of minimizers]\label{prop:existence}
For every $\rho,\eta\in\cM_+(\cX)$, the infimum in \eqref{eq:static-wfr} is attained. If $\rho=\eta=0$, the unique effective minimizer is $a=b=0$. If one endpoint is zero, say $\eta=0$, the value is
\begin{equation}\label{eq:wfr-zero}
    \WFR_\delta^2(\rho,0)=4\delta^2|\rho|.
\end{equation}
\end{proposition}
\begin{proof}
The proof is provided in Supplementary Material~\ref{app:proof-existence}.   
\end{proof}

\begin{proposition}[Elementary properties]\label{prop:basic-properties}
For all $\rho,\eta\in\cM_+(\cX)$ and $\lambda\geq0$, the discrete WFR functional satisfies:
\begin{align}
    &\WFR_\delta(\rho,\eta)\geq0,\label{eq:nonnegative}\\
    &\WFR_\delta(\rho,\eta)=\WFR_\delta(\eta,\rho),\label{eq:symmetric}\\
    &\WFR_\delta^2(\lambda\rho,\lambda\eta)=\lambda \WFR_\delta^2(\rho,\eta),\label{eq:homogeneous}\\
    &2\delta\left|\sqrt{|\rho|}-\sqrt{|\eta|}\right|
      \leq \WFR_\delta(\rho,\eta)
      \leq 2\delta\sqrt{|\rho|+|\eta|}.\label{eq:mass-bounds}
\end{align}
Moreover, $\WFR_\delta(\rho,\eta)=0$ if and only if $\rho=\eta$.
\end{proposition}
\begin{proof}
The proof is provided in Supplementary Material~\ref{app:proof-basic-properties}.   
\end{proof}

\begin{theorem}[Metric property]\label{thm:metric}
When \(\max_{i,j}C_{ij} < \pi^2\delta^2\), and \(C_{ij}^{\frac{1}{p}}\) is a metric on \(\cX\), \((1\le p\le2)\), the map
\begin{equation*}
    (\rho,\eta)\mapsto \WFR_\delta(\rho,\eta)
\end{equation*}
is a metric on $\cM_+(\cX)$.
\end{theorem}
\begin{proof}
The proof is provided in Supplementary Material~\ref{app:proof-wfr-metric}.    
\end{proof}

Similar to the continuous space case, we can also derive an equivalent optimal entropy transport formulation for discrete WFR problem \cite{liero2018optimal,chizat2018unbalanced,wfrfm2025}. Such a form is easier to solve numerically.
\begin{proposition}[Optimal entropy transport form]\label{prop:entropy}
Let
\begin{equation}\label{eq:et-cost}
    c_{ij}=
    -2\log\left[\cos\left(\dfrac{\sqrt{C_{ij}}}{2\delta}\right)\right]
\end{equation}
For nonnegative vectors $r,s$, define the KL-divergence
\begin{equation}\label{eq:kl}
    \KL(r\,|\,s)=
    \sum_i
    \left[
        r_i\log\left(\frac{r_i}{s_i}\right)-r_i+s_i
    \right],
\end{equation}
with the usual conventions $0\log0=0$ and $\KL(r\,|\,s)=+\infty$ if $r_i>0$ and $s_i=0$ for some $i$. When \(\max_{i,j}C_{ij} < \pi^2\delta^2\)
\begin{equation}\label{eq:et-wfr}
    \frac{1}{4\delta^2}\WFR_{\delta}^2(\rho,\eta)
    =
    \inf_{\gamma\in\Rnum_+^{N\times N}}
    \left\{
        \sum_{i,j=1}^N c_{ij}\gamma_{ij}
        +
        \KL(\gamma\one\,|\,\rho)
        +
        \KL(\gamma^T\one\,|\,\eta)
    \right\}.
\end{equation}
\end{proposition}
\begin{proof}
The proof is provided in Supplementary Material~\ref{app:proof-entropy}.   
\end{proof}

\begin{proposition}[Semi-coupling]\label{prop:et-kkt}
Assume $\rho_i,\eta_j\ge0$ and $|\rho|,|\eta|>0$. Let $\gamma^\ast$ solve \eqref{eq:et-wfr}, and set \(r_i^\ast=(\gamma^\ast\one)_i,\quad s_j^\ast=((\gamma^\ast)^T\one)_j.\) Then the optimal semi-coupling in Definition~\ref{def:staticwfr} can be recovered from $\gamma^\ast$ by
\begin{equation}\label{eq:recover-ab}
    a_{ij}^\ast=\gamma_{ij}^\ast\frac{\rho_i}{r_i^\ast},
    b_{ij}^\ast=\gamma_{ij}^\ast\frac{\eta_j}{s_j^\ast},\  
    \rho_i\eta_j\neq0,\qquad 
    a_{ij}^\ast = b_{ij}^\ast = \gamma_{ij}^\ast = 0,\ 
    \rho_i\eta_j=0.
\end{equation}
\end{proposition}
\begin{proof}
The proof is provided in Supplementary Material~\ref{app:proof-et-kkt}.   
\end{proof}

The optimal entropy transport objective in \eqref{eq:et-wfr} is convex in $\gamma$, and it can be efficiently solved by Sinkhorn-like iteration \cite{chizat2018unbalanced}. With proposition \ref{prop:et-kkt}, we can solve the static WFR semi-coupling \((a,b)\) efficiently.

\subsection{From static to dynamic}
Given the optimal semi-coupling \((a,b)\) between two nonnegative measures \(\rho,\eta\in\mathcal{M}_+(\mathcal{X})\), we can derive the travelling Diracs \(m_{ij}(t),k_{ij}(t)\) from the boundary condition \(m_{ij}(0)=a_{ij},m_{ij}(1)=b_{ij}\). The dynamic measure flow corresponding to the dynamic WFR problem is then defined.
\begin{equation}
\label{eq:dynamic discrete WFR flow}
    \mu_t
    =
    \sum_{i,j}m_{ij}(t)
    \left[
        (1-k_{ij}(t))\delta_{x_i}
        +
        k_{ij}(t)\delta_{x_j}
    \right].
\end{equation}

Intuitively, the dynamic WFR flow is the average of travelling Diracs over semi-coupling. In continuous scenarios, a standard algorithm for realizing this average is flow matching \cite{lipman2023flow}. Previous works have utilized flow matching and its unbalanced extension to solve dynamic OT, WFR, and other unbalanced OT extensions \cite{cfm_tong,wang2025joint,wfrfm2025,USB,SUDO} in Euclidean spaces. Although flow matching is developed for continuous spaces, recent works have already proposed discrete flow matching for learning dynamics in discrete spaces \cite{campbell2024generative,gat2024discrete}. 

Another advantage of flow matching is that it does not directly learn the measure flow \(\mu_t\), but learn the velocity field \(v_t\) and growth rate \(g_t\) which generate \(\mu_t\), i.e. the dynamics. In discrete scenarios, the dynamics can be represented by a state dependent growth rate and a transition kernel. Unlike continuous dynamic WFR, the dynamic discrete WFR defined above is not a Markov process. But we can at least learn its Markovian reconstruction, which is a Markov process that generates \(\mu_t\). We discuss the details in the next section.

\section{Markovian reconstruction of WFR and unbalanced discrete flow matching}
\label{sec:unbalanced_dfm}
In this section, we introduce a general framework for unbalanced discrete flow matching on a finite state space $\mathcal X$. Based on it, we derived the closed form of the Markovian reconstruction of the discrete WFR dynamics. In contrast to standard discrete flow matching \cite{gat2024discrete}, where the total mass of the evolving distribution is preserved, the present formulation allows both transport across discrete states and state-dependent mass growth or decay. This provides a natural modeling framework for discrete dynamics with birth and death. We first introduce the marginalization theory, and developed the closed form of the dynamic discrete WFR flow. For large \(|\mathcal{X}|\), we also developed an unbalanced discrete conditional flow matching objective for learning the solution which may be intractable.

\subsection{Unbalanced Discrete Dynamics}
\label{subsec:unbalanced_discrete_dynamics}
We first describe the stochastic process underlying the unbalanced discrete dynamics considered in this work. Let $\mathcal X$ be a finite state space, and let $(x_t)_{t\in[0,1]}$ be a time-inhomogeneous continuous-time Markov chain on $\mathcal X$. For $x,y\in\mathcal X$ with $y\neq x$, we denote by $ q_t(x,y)\geq 0 $ the transition rate from state $x$ to state $y$ at time $t$.
To model
non-conservative dynamics, we further attach to each trajectory a nonnegative mass
process $(m_t)_{t\in[0,1]}$ satisfying
\begin{equation}
\frac{d}{dt}\log m_t = g_t(x_t)
\label{eq:trajectory_mass_growth}
\end{equation}
where $g_t(x)$ is a state-dependent growth rate: positive values increase the mass carried by trajectories at state $x$, while negative values decrease it.

For each state $x\in\mathcal X$, we define
\begin{equation}
p_t(x)
=
\mathbb E\left[m_t\mathbf 1_{\{x_t=x\}}\right].
\label{eq:weighted_law_definition}
\end{equation}
This is the expected mass carried by trajectories that are located at state
$x$ at time $t$. Hence $p_t(x)$ plays the role of the state-wise mass density in
the discrete space, or say the nonnegative measure on \(\mathcal{X}\). In contrast to a probability distribution, $p_t$ is not
required to sum to one:
\[
\sum_{x\in\mathcal X}p_t(x)
=
\mathbb E[m_t].
\]
The total mass may increase or decrease over time, depending on the
state-dependent growth rate $g_t$.

Under the above assumptions, the weighted marginal satisfies the component-wise
unbalanced Kolmogorov forward equation
\begin{equation}
\partial_t p_t(x)
=
\sum_{y\neq x}p_t(y)q_t(y,x)
-
p_t(x)\sum_{y\neq x}q_t(x,y)
+
p_t(x)g_t(x),
\qquad x\in\mathcal X.
\label{eq:unbalanced_kfe}
\end{equation}
The first term describes the incoming mass flux into state $x$, the second term
describes the outgoing mass flux from state $x$, and the last term describes
local multiplication or decay of mass at state $x$. Hence the transition rates
$q_t(x,y)$ determine how mass is transported across states, whereas the growth
rate $g_t(x)$ determines how mass is created or removed locally.

\subsection{Conditional and Marginal Dynamics}
\label{subsec:conditional_marginal_dynamics}
We now describe how travelling Diracs can be marginalized into the marginal dynamic flow. In the unbalanced setting, a conditional measure flow is specified by both endpoint states and endpoint masses. We write
\[
z=((x_0,m_0),(x_1,m_1))
\in
\mathcal Z:=(\mathcal X\times\mathbb R_{\ge0})^2,
\]
where $x_0,x_1\in\mathcal X$ and $m_0,m_1\ge0$. Let $\pi$ be the distribution on $\mathcal Z$. For each $z\in\mathcal Z$, we construct a conditional flow by separating the
total mass from the normalized discrete distribution. Let
\[
m_t(z)\ge0,
\qquad
\widetilde p_t(\cdot\mid z)\in\mathcal P(\mathcal X),
\]
with boundary conditions
\[
m_0(z)=m_0,\quad m_1(z)=m_1,
\qquad
\widetilde p_0(x\mid z)=\mathbf 1_{\{x=x_0\}},
\quad
\widetilde p_1(x\mid z)=\mathbf 1_{\{x=x_1\}}.
\]
The normalized path evolves conservatively according to
\begin{equation}
\partial_t \widetilde p_t(x\mid z)
=
\sum_{y\neq x}
\widetilde p_t(y\mid z)q_t(y,x\mid z)
-
\widetilde p_t(x\mid z)
\sum_{y\neq x}q_t(x,y\mid z),
\qquad x\in\mathcal X.
\label{eq:udfm_normalized_conditional_forward}
\end{equation}
We then define the unbalanced conditional flow by
\begin{equation}
p_t(x\mid z)
=
m_t(z)\widetilde p_t(x\mid z).
\label{eq:udfm_conditional_path_construction}
\end{equation}
Thus
\[
p_0(x\mid z)=m_0\mathbf 1_{\{x=x_0\}},
\qquad
p_1(x\mid z)=m_1\mathbf 1_{\{x=x_1\}}.
\]

The conditional growth rate is assigned to the mass dynamics:
\begin{equation}
g_t(x\mid z)
=
g_t(z)
:=
\partial_t\log m_t(z),
\qquad x\in\mathcal X.
\label{eq:udfm_conditional_growth_scalar}
\end{equation}
Differentiating \eqref{eq:udfm_conditional_path_construction} and using
\eqref{eq:udfm_normalized_conditional_forward} yields
\begin{equation}
\partial_t p_t(x\mid z)
=
\sum_{y\neq x}p_t(y\mid z)q_t(y,x\mid z)
-
p_t(x\mid z)\sum_{y\neq x}q_t(x,y\mid z)
+
p_t(x\mid z)g_t(z).
\label{eq:conditional_ukfe}
\end{equation}
The marginal flow is obtained by averaging over the endpoint coupling:
\begin{equation}
p_t(x)
=
\int_{\mathcal Z}p_t(x\mid z)\,\pi(z)\dd z,
\qquad x\in\mathcal X.
\label{eq:udfm_marginal_path}
\end{equation}
For every $x$ with $p_t(x)>0$, define the marginal transition and growth rates by
\begin{align}
q_t(x,y)
&=
\frac{
\int_{\mathcal Z} q_t(x,y\mid z)p_t(x\mid z)\,\pi(z)\dd z
}{
p_t(x)
},
\qquad y\neq x,
\label{eq:udfm_marginal_rate}
\\
g_t(x)
&=
\frac{
\int_{\mathcal Z} g_t(z)p_t(x\mid z)\,\pi(z)\dd z
}{
p_t(x)
}.
\label{eq:udfm_marginal_growth}
\end{align}
When $p_t(x)=0$, the values of $q_t(x,y)$ and $g_t(x)$ can be chosen
arbitrarily.

Under the above definitions, the marginal flow $p_t(x)$ satisfies
\eqref{eq:unbalanced_kfe}, as follows by differentiating
\eqref{eq:udfm_marginal_path} and applying
\eqref{eq:conditional_ukfe} together with
\eqref{eq:udfm_marginal_rate}--\eqref{eq:udfm_marginal_growth}. Thus the
marginal dynamics is induced by averaging conditional dynamics, with effective transition and growth rates weighted by the conditional mass at the current state. 

\subsection{The closed form of dynamic discrete WFR}
Given the optimal semi-coupling of static WFR (\ref{eq:static-wfr}) \((a,b)\), we can sample the conditional variable \(z\) with probability \(\text{P}\Big(((x_i,a_{ij}),(x_j,b_{ij}))\Big)\propto a_{ij}+b_{ij},\forall i,j\), and use the corresponding travelling Dirac with boundary mass \((|a|+|b|)(\frac{a_{ij}}{a_{ij}+b_{ij}},\frac{b_{ij}}{a_{ij}+b_{ij}})\) as the conditional flow \(p_t(x|z)\). Then, the marginal measure path (\ref{eq:udfm_marginal_path}) is exactly the dynamic WFR flow (\ref{eq:dynamic discrete WFR flow}). In this perspective, we can view the marginal transition and growth rate \(q_t,g_t\) as the Markovian reconstruction dynamics of the dynamic WFR flow.

In this case, the conditional variable \(z\) can be viewed as a discrete random variable with distribution \(\pi(z)\). Therefore, the marginalization relation (\ref{eq:udfm_marginal_path},\ref{eq:udfm_marginal_rate},\ref{eq:udfm_marginal_growth}) are actually not intractable integrals, but summations.

\begin{equation}
\label{eq: analytic solution}
\begin{aligned}
p_t(x)
&=
\sum_{z\in\mathcal Z}p_t(x\mid z)\,\pi(z),
\qquad x\in\mathcal X,\\
q_t(x,y)
&=
\frac{
\sum_{z\in\mathcal Z} q_t(x,y\mid z)p_t(x\mid z)\,\pi(z)
}{
p_t(x)
},\\
g_t(x)
&=
\frac{
\sum_{z\in\mathcal Z} g_t(z)p_t(x\mid z)\,\pi(z)
}{
p_t(x)
}.
\end{aligned}
\end{equation}

With \(\pi\) determined, it is sufficient to derive the conditional \(q_t,g_t\) from the travelling Diracs. When \(z=((x_i,a_{ij}),(x_j,b_{ij}))\), the conditional flow is
\begin{equation}
p_t(\cdot|z)=m_t(z)((1-k_t(z))\delta_{x_i}+k_t(z)\delta_{x_j})
\end{equation}
where \(m_t(z),k_t(z)\) are determined by (\ref{eq:travelling Dirac}) with boundary \((|a|+|b|)(\frac{a_{ij}}{a_{ij}+b_{ij}},\frac{b_{ij}}{a_{ij}+b_{ij}})\). Using (\ref{eq:udfm_normalized_conditional_forward},\ref{eq:udfm_conditional_growth_scalar}), it is sufficient to set \(g_t(z)=\partial_t\ln m_t(z)\), and find a conditional transition rates \(q_t(x,y|z)\) such that
\begin{equation}
\partial_t \widetilde p_t(x\mid z)
=
\sum_{y\neq x}
\widetilde p_t(y\mid z)q_t(y,x\mid z)
-
\widetilde p_t(x\mid z)
\sum_{y\neq x}q_t(x,y\mid z)
\end{equation}
where \(\widetilde p_t(x\mid z)=(1-k_t(z))\delta_{x_i}+k_t(z)\delta_{x_j}\). By direct calculation, we have
\begin{equation}
q_t(x_i,x_j|z)=\frac{\partial_tk_t(z)}{1-k_t(z)},i\neq j,\quad q_t(x_i,y|z)=0,y\neq x_i, x_j, \quad q_t(x,y|z)=0,x\neq x_i
\end{equation}
Since \(q_t\) is a Markov generator, the diagonal elements are automatically determined by the off-diagonal elements.

Therefore, when \(N=|\mathcal{X}|\) is not so  large, the summations are tractable. We can directly construct the closed form of discrete WFR via (\ref{eq: analytic solution}).

\subsection{Unbalanced discrete conditional flow matching objectives}
\label{subsec:udfm_losses}
When \(N=|\mathcal{X}|\) is too large to calculate the sums (\ref{eq: analytic solution}), or for more general cases beyond WFR where \(z\) may be a continuous random variable, the marginals (\ref{eq:udfm_marginal_path},\ref{eq:udfm_marginal_rate},\ref{eq:udfm_marginal_growth}) are intractable. In order to approximate them numerically, we developed an unbalanced conditional discrete flow matching framework. Like conventional flow matching algorithms, one can implicitly regress the marginals by regressing the tractable conditionals.

We now define the learning objectives for unbalanced conditional discrete flow matching.
The goal is to learn marginal transition rates and growth rates that generate the
marginal target path. We parameterize the off-diagonal transition rates by
$q_t^\theta(x,y)$, $y\neq x$, and the growth rate by $g_t^\phi(x)$. Since the
diagonal entries of a Markov generator are determined by the off-diagonal
rates, we only match the outgoing off-diagonal rate vector. For
$x\in\mathcal X$, define
\begin{equation}
\langle u,v\rangle_x
:=
\sum_{y\neq x}u(y)v(y),
\qquad
\|v\|_{2,x}^2
:=
\langle v,v\rangle_x .
\label{eq:udfm_offdiag_inner_product}
\end{equation}
Here $q_t(x,\cdot)$ denotes the vector
$\{q_t(x,y):y\neq x\}$.

\paragraph{Marginal UDFM objective}
The marginal unbalanced discrete flow matching objective is a direct regression loss
\begin{equation}
\begin{aligned}
\mathcal L_{\mathrm{UDFM}}(\theta,\phi)
&=
\mathbb E_{t\sim U[0,1]}
\sum_{x\in\mathcal X}p_t(x)
\Bigl[(1-t)^2
\big\|
q_t^\theta(x,\cdot)-q_t(x,\cdot)
\big\|_{2,x}^2\\
&+
\lambda
\big|
g_t^\phi(x)-g_t(x)
\big|^2
\Bigr],
\label{eq:udfm_loss}
\end{aligned}
\end{equation}
where $q_t(x,y)$ and $g_t(x)$ are the marginal target fields defined in
\eqref{eq:udfm_marginal_rate} and \eqref{eq:udfm_marginal_growth}.

\paragraph{Conditional UDFM objective.}
In practice, the marginal measure path and target fields may be intractable due to the integral. Instead, one can train directly against conditional bridge targets. This gives
\begin{equation}
\begin{aligned}
\mathcal L_{\mathrm{CUDFM}}(\theta,\phi)
&=
\mathbb E_{t\sim U[0,1]}
\int_{\mathcal Z}
\sum_{x\in\mathcal X}p_t(x\mid z)
\Bigl[(1-t)^2
\big\|
q_t^\theta(x,\cdot)-q_t(x,\cdot\mid z)
\big\|_{2,x}^2
\\
&\hspace{2.3cm}
+
\lambda
\big|
g_t^\phi(x)-g_t(x\mid z)
\big|^2
\Bigr]
\pi(z)\dd z.
\end{aligned}
\label{eq:cudfm_loss}
\end{equation}
The next result shows that this conditional objective has the same gradient as the marginal objective.

\begin{theorem}[Equivalence of UDFM and CUDFM]
\label{thm:udfm_cudfm_equivalence}
Assume that $p_t(x)>0$ for all $x\in\mathcal X$ and $t\in[0,1]$. Let the
marginal target fields be defined by
\eqref{eq:udfm_marginal_rate} and \eqref{eq:udfm_marginal_growth}. Then
\begin{equation}
\mathcal L_{\mathrm{UDFM}}(\theta,\phi)
=
\mathcal L_{\mathrm{CUDFM}}(\theta,\phi)
+
C,
\label{eq:udfm_cudfm_equivalence}
\end{equation}
where $C$ is independent of $(\theta,\phi)$. Consequently,
\begin{equation}
\nabla_{\theta,\phi}\mathcal L_{\mathrm{UDFM}}(\theta,\phi)
=
\nabla_{\theta,\phi}\mathcal L_{\mathrm{CUDFM}}(\theta,\phi).
\label{eq:udfm_cudfm_gradient_equivalence}
\end{equation}
\end{theorem}
\begin{proof}
The proof is provided in Supplementary Material~\ref{app:proof-udfm-equivalence}.
\end{proof}

\begin{remark}
We did not use the cross-entropy loss used by the original discrete flow matching \cite{gat2024discrete} here, because the loss equivalence for it only holds when all conditional flows share the same progress variable \(k(t)\), which is not the case for WFR. For discrete OT where all travelling Diracs share the same \(k(t)=t\), one can use the original discrete flow matching.
\end{remark} 
\begin{remark}
In practice, \(q_t(x,\cdot\mid z) = O(\frac{1}{1-t})\) as \(t\to 1\). The weight \((1-t)^2\) is therefore added to avoid numerical instability.
\end{remark} 

The equivalence theorem justifies that training with conditional loss is equivalent to training with the marginal one. Therefore, one can regress the intractable marginal targets by regressing the tractable conditional targets. Given the dynamics \(q_t(x,y),g_t(x)\), one can recover the measure path \(p_t(x)\) by simulating the master equation (\ref{eq:unbalanced_kfe}) from initial value \(p_0(x)\).

\section{Multiscale consistency}
\label{sec:coarse_consistency}
In practice, cell state annotations are often available at different resolutions. Therefore, a natural question is whether the reconstructed dynamics is consistent across multiscale annotations, or say, under state aggregation. Let $\mathcal X=\{1,\ldots,N\}$ be a fine state space and $\bar{\mathcal X}=\{1,\ldots,M\}$ a coarse state space. A surjective map \(P:\mathcal X\to\bar{\mathcal X}\) induces the partition $B_A=P^{-1}(A)$ and the pushforward \((P_\#\rho)_A=\sum_{i\in B_A}\rho_i .\) Let $\rho,\eta\in\mathcal M_+(\mathcal X)$, and write
$\bar\rho=P_\#\rho$ and $\bar\eta=P_\#\eta$. For convenience, define \(
\zeta_{ij}
=
\cos\!\left(\frac{\sqrt{C_{ij}}}{2\delta}\right),
\bar\zeta_{AB}
=
\cos\!\left(\frac{\sqrt{\bar C_{AB}}}{2\delta}\right),
\) and assume $C_{ij},\bar{C}_{AB}<\delta^2\pi^2$.
The static WFR objective on the fine space can then be written as
\[
\mathcal J_\zeta(a,b)
=
\sum_{i,j}
\left(
a_{ij}+b_{ij}
-2\zeta_{ij}\sqrt{a_{ij}b_{ij}}
\right),
\]
under
\[
\sum_j a_{ij}=\rho_i,
\qquad
\sum_i b_{ij}=\eta_j.
\]
The squared WFR distance equals
$4\delta^2\min_{a,b}\mathcal J_\zeta(a,b)$.

Given a fine semi-coupling $(a,b)$, define its aggregation by
\begin{equation}
\bar a_{AB}
=
\sum_{i\in B_A,j\in B_B}a_{ij},
\qquad
\bar b_{AB}
=
\sum_{i\in B_A,j\in B_B}b_{ij}.   
\end{equation}

\begin{theorem}[Exact consistency under block-homogeneous costs]
\label{thm:exact_coarse}
Assume
\begin{equation}
\label{eq:cost consistency}
C_{ij}
=
\bar C_{P(i)P(j)}
\qquad \forall i,j.    
\end{equation}

Then
\begin{equation}
\mathrm{WFR}_{\delta}(\rho,\eta)
=
\mathrm{WFR}_{\delta}(\bar\rho,\bar\eta).  
\end{equation}
Moreover, if $(a,b)$ is any fine-scale optimal semi-coupling, its aggregation $(\bar a,\bar b)$ is coarse-scale optimal, and the
corresponding dynamic WFR flow satisfy
\begin{equation}
\label{eq:lift}
P_\#\mu_t^{\mathrm{fine}}
=
\mu_t^{\mathrm{coarse}},
\qquad t\in[0,1],    
\end{equation}
where the coarse path is reconstructed from $(\bar a,\bar b)$.
Conversely, every coarse optimal semi-coupling admits a fine optimal
lifting satisfying \eqref{eq:lift}.
\end{theorem}

\begin{proof}
For any fine feasible $(a,b)$, its aggregation is coarse feasible.
For each block pair $(A,B)$, Cauchy--Schwarz gives
\[
\sum_{i\in B_A,j\in B_B}\sqrt{a_{ij}b_{ij}}
\leq
\sqrt{\bar a_{AB}\bar b_{AB}}.
\]
Using \eqref{eq:cost consistency},
\[
\bar{\mathcal J}_{\bar\zeta}(\bar a,\bar b)
\leq
\mathcal J_\zeta(a,b),
\]
and hence the coarse optimal value is no larger than the fine one.

For the other direction, let $(\bar a,\bar b)$ be coarse feasible.
For $i\in B_A$ and $j\in B_B$, define
\[
p_i^A=\frac{\rho_i}{\bar\rho_A},
\qquad
q_j^B=\frac{\eta_j}{\bar\eta_B},
\]
whenever the denominators are positive; if a block has zero total mass, choose an arbitrary probability vector on that block. Set
\begin{equation}
\label{eq:fine_to_coarse}
a_{ij}=p_i^Aq_j^B\,\bar a_{AB},
\qquad
b_{ij}=p_i^Aq_j^B\,\bar b_{AB}.   
\end{equation}

This is a fine feasible semi-coupling, and we have
\[
\sqrt{a_{ij}b_{ij}}
=
p_i^Aq_j^B\sqrt{\bar a_{AB}\bar b_{AB}}.
\]
Summing over each block pair \((A,B)\) gives
\[
\mathcal J_\zeta(a,b)=\bar{\mathcal J}_{\bar\zeta}(\bar a,\bar b).
\]
Thus, \(\mathrm{WFR}_{\delta}(\rho,\eta)
=
\mathrm{WFR}_{\delta}(\bar\rho,\bar\eta).\)

Now let $(a,b)$ be fine optimal. Since its aggregation is also optimal,
all blockwise Cauchy--Schwarz inequalities above must be equalities.
Hence, for every active block pair $(A,B)$, there exist
$\lambda_{ij}\geq0$ with $\sum_{i,j}\lambda_{ij}=1$ such that
\begin{equation}
a_{ij}=\lambda_{ij}\bar a_{AB},
\qquad
b_{ij}=\lambda_{ij}\bar b_{AB}.
\end{equation}

Since the travelling Dirac is positively homogeneous in its endpoint masses, we have
\begin{equation}
\mu_t^{ij}=m_{ij}(t)[(1-k_{ij}(t))\delta_{x_i}+k_{ij}(t)\delta_{x_j}]=\lambda_{ij}m_{AB}(t)[(1-k_{AB}(t))\delta_{x_i}+k_{AB}(t)\delta_{x_j}]
\end{equation}
Therefore,
\[
\sum_{i\in B_A,j\in B_B}
P_\#\mu_t^{ij}
=\sum_{i\in B_A,j\in B_B}\lambda_{ij}\cdot\Big(m_{AB}(t)[(1-k_{AB}(t))\delta_A+k_{AB}(t)\delta_B]\Big)=
\mu_t^{AB},
\]
and summing over all block pairs yields \eqref{eq:lift}.
The lifting \eqref{eq:fine_to_coarse} gives the converse.
\end{proof}
This theorem shows that, under the block-homogeneous cost condition \eqref{eq:cost consistency} , the coarse-scale dynamics and the fine-scale dynamics correspond exactly. However, in practice, fine states within one coarse class generally have slightly different pairwise costs. We next quantify the resulting error. We define the geometric distortion
\begin{equation}
\label{eq:geometric distortion}
\varepsilon
=
\max_{i,j}
\left|
C_{ij}
-
\bar C_{P(i)P(j)}
\right|.    
\end{equation}
By direct calculation, we have
\begin{equation}
\begin{aligned}
\left|\zeta_{ij}-\bar\zeta_{P(i)P(j)}\right|&=\left|\cos \frac{\sqrt{C_{ij}}}{2\delta}-\cos\frac{\sqrt{\bar C_{P(i)P(j)}}}{2\delta}\right|\\&\le\frac{1}{2}\left|\big(\frac{\sqrt{C_{ij}}}{2\delta}\big)^2-\big(\frac{\sqrt{\bar C_{P(i)P(j)}}}{2\delta}\big)^2\right|=\frac{1}{8\delta^2}\left|C_{ij}-\bar C_{P(i)P(j)}\right|\le\frac{\varepsilon}{8\delta^2}
\end{aligned}
\end{equation}

\begin{proposition}[Total mass stability under approximate aggregation]
\label{prop:approx_coarse}
Let \(R=|\rho|,S=|\eta|.\) Then
\begin{equation}
\label{eq:wfr cost consistency}
\left|
\mathrm{WFR}_{\delta}^2(\rho,\eta)
-
\mathrm{WFR}_{\delta}^2(\bar\rho,\bar\eta)
\right|
\leq\varepsilon\sqrt{RS}.    
\end{equation}

Furthermore, if $M_f(t)$ and $M_c(t)$ denote the total masses along
any optimal fine and coarse WFR flows, respectively, then
\begin{equation}
\label{eq:mass consistency}
|M_f(t)-M_c(t)|
\leq
\frac{t(1-t)}{4\delta^2}\varepsilon\sqrt{RS}    
\end{equation}
and, for $R,S>0$,
\begin{equation}
\label{eq:RME consistency}
\sup_{t\in[0,1]}
\frac{|M_f(t)-M_c(t)|}{M_c(t)}
\leq
\frac{\varepsilon}{8\delta^2}.   
\end{equation}
\end{proposition}
\begin{proof}
The proof is provided in Supplementary Material~\ref{app:proof-total-mass-stability}.
\end{proof}

Proposition~\ref{prop:approx_coarse} provides a quantitative characterization of the relative mass error (RME) between the coarse-scale and fine-scale dynamics, showing that the error is bounded linearly in terms of the geometric distortion \eqref{eq:geometric distortion}. For the stability of the full state-composition measure flow, we have the following convergence proposition.

\begin{proposition}[Measure flow convergence under\allowbreak{} state aggregation]
\label{prop:path_stability}
Let $\{C^{(n)}\}_{n\geq1}$ be a sequence of fine-scale costs satisfying \(C_{ij}^{(n)},\bar C_{AB}<\delta^2\pi^2\) and a sequence of geometric distortion
\begin{equation}
\label{eq:epsilon_n}
\varepsilon_n
:=
\max_{i,j}
\left|
C_{ij}^{(n)}
-
\bar C_{P(i)P(j)}
\right|
\longrightarrow0.
\end{equation}
For each $n$, let $(a^{(n)},b^{(n)})$ be an arbitrary optimal
fine-scale semi-coupling under the cost $C^{(n)}$, and let
$\mu_t^{(n)}$ be its corresponding dynamic WFR flow. Then every subsequential limit of $P_\#\mu_t^{(n)}$, with respect to
the norm \(\sup_{t\in[0,1]}\|\cdot\|_1\), is an optimal coarse-scale WFR flow. In particular, if the coarse optimal WFR flow
$\bar\mu_t^*$ is unique, then
\begin{equation}
\label{eq:path consistency}
\sup_{t\in[0,1]}
\left\|
P_\#\mu_t^{(n)}
-
\bar\mu_t^*
\right\|_1
\longrightarrow0.
\end{equation}
If $R,S>0$, then the normalized measure flows also converge uniformly
in total variation:
\begin{equation}
\label{eq:TV consistency}
\sup_{t\in[0,1]}
\mathrm{TV}
\left(
\frac{P_\#\mu_t^{(n)}}{|P_\#\mu_t^{(n)}|},
\frac{\bar\mu_t^*}{|\bar\mu_t^*|}
\right)
\longrightarrow0.
\end{equation}
\end{proposition}
\begin{proof}
The proof is provided in Supplementary Material~\ref{app:proof-measure-flow-convergence}.
\end{proof}

These results show that discrete WFR commutes exactly with state aggregation under block-homogeneous condition and remains stable under small within-block geometric distortion. This provides a theoretical justification for comparing fine-to-coarse inference with direct coarse inference in Sec.~\ref{sec:consistency}. In practice, discrete WFR modeling is theoretically robust to cell type annotation resolution.

\section{Numerical results}
\label{sec:experiments}
In this section, we compute dynamic discrete WFR flows on several synthetic and real world datasets. Unless otherwise specified, we use the closed form solution directly and perform inference by integrating the master equation. The corresponding algorithm is called WFR-DFM (WFR discrete flow matching).

For each biological dataset, we are given several cell by gene transcriptomic matrices, each corresponding to an observation time point. Each matrix entry records the RNA count of a particular gene in a particular cell. We first reduce the dimensionality of the data and then group the cells into a finite set of discrete states according to the provided cell type annotations. The position of each state in the reduced dimensional space is defined as the centroid of its constituent cells. The edge weight between two states is defined as the squared Euclidean distance between their centroids. Finally, the measure at each observation time point is naturally determined by the total number of observed cells and their proportions across the discrete states. The details of synthetic data generation and real data processing are left to the supplementary material.

\subsection{The necessity of WFR}
In biological processes, cells not only shift their gene expression, but also divide or die. To faithfully recover the cellular dynamics, it is necessary to simultaneously model both state transition and population growth. To demonstrate the necessity of applying WFR rather than balanced OT, we calculate the dynamics induced by OT and WFR on a four-state synthetic system with states $A,B,C,D$. The ground-truth
measure was generated by
\begin{equation}
\label{eq: toy dynamics}
  \dot{\mu}_t
  =
  Q_{\mathrm{true}}^\top\mu_t
  +
  g_{\mathrm{true}}\odot\mu_t,
\end{equation}
with
\begin{equation}
\begin{aligned}
Q_{\mathrm{true}}
&=
\begin{pmatrix}
-0.6 & 0.6 & 0 & 0\\
0 & 0 & 0 & 0\\
0 & 0 & 0 & 0\\
0 & 0 & 0 & 0
\end{pmatrix},\\
g_{\mathrm{true}}&=(0,0,-1.4,1.6)^\top,
\qquad
\mu_0=(0.45,0.15,0.38,0.02)^\top .
\end{aligned}
\end{equation}
That is, the only transport edge is $A\to B$ with transport rate \(0.6\). State $C$ is isolated
from transport and undergoes pure death, whereas state $D$ is isolated from
transport and undergoes pure growth. The observations \(\mu_0,\mu_1,\mu_2\) were taken at $t=0,1,2$. The resulting total masses $1.0000$, $0.7928$, and $1.1138$. The ground cost was defined by $C_{ii}=0$, $C_{AB}=C_{BA}=0.25$, and
$C_{ij}=5$ for all remaining $i\neq j$. Ground-truth transition flows on edge \(i\to j\) were obtained by integrating $\mu_i(t)q_{ij}(t)$.

We calculate the WFR dynamics \((\delta=1)\) and the balanced OT dynamics in that synthetic system, and predict the measures at the time point \(1,2\), denoted as \(\hat{\mu}_1,\hat{\mu}_2\). The predictions were obtained by using \(\mu_0\) as initial measure, and integrating (\ref{eq: toy dynamics}) with predicted \(Q_t,g_t\). To measure the accuracy of measure prediction, we calculated the total variation (TV) between normalized true measure and normalized predicted measure, and the relative mass error (RME) between the two measures.
\begin{equation}
\text{TV}_{\text{normalized}}(\rho,\eta)=\frac{1}{2}\sum_{x\in\mathcal{X}}\Big|\frac{\rho(x)}{|\rho|}-\frac{\eta(x)}{|\eta|}\Big|,\qquad \text{RME}(\rho,\eta) = \frac{||\rho|-|\eta||}{|\eta|}.
\end{equation}
TV quantifies the difference in shape between the two measures, whereas RME quantifies the difference in their total mass. As shown in Table \ref{tab:toy system}, macroscopically, both OT and WFR can recover the normalized measure, while OT fails to recover the total mass since it did not model the growth. The microscopic transport comparison provides evidence independent of total mass prediction.
Balanced OT assigned $83\%$ of its interval flows,
respectively, to edges outside the true $A\to B$ support, thereby representing
the death of $C$ and growth of $D$ as spurious transport. In contrast, WFR
recovered only the $A\to B$ edge for every tested value, with no spurious flow ($10^{-11}$). This microscopic result shows that although both methods can connect the snapshots accurately, they achieve this by different dynamics. Balanced OT tends to compensate the birth-death dynamics by spurious transports. Therefore, it is necessary to simultaneously model the growth rate, in order to recover the true transport dynamics. 

\begin{table}[h]
\centering
\caption{Comparison of WFR and balanced OT on synthetic system. TV, RME are averaged over timepoint 1, 2.}
\label{tab:toy system}
\begin{tabular}{lcc}
\toprule
Metric & WFR & Balanced OT\\
\midrule
TV $\downarrow$
  & 0.0211 & 0.0304\\
RME $\downarrow$
  & 0.0273 & 0.1614\\
False edge flow fraction $\downarrow$
  & $8.04\times10^{-11}$ & 0.8270\\
\bottomrule
\end{tabular}
\end{table}

\subsection{Unbalanced discrete flow matching recovers the solution}

We evaluated whether unbalanced discrete conditional flow matching (CFM) recovers the discrete WFR measure path (Direct) on a seven-state real EMCC data, a 50-state synthetic problem, a 1,000-state synthetic problem. Each dataset contains four observation time points, and each interval between two consecutive time points is uniformly divided into ten subintervals. Excluding the initial time point, this yields 30 evaluation time points in total. We train CFM using the four observed time points, then perform inference from time zero following the learned dynamics and evaluate the resulting trajectories at these 30 time points. CFM was trained for 1,000 steps on the 50-state problem and EMCC and for 10,000 steps on the 1,000-state problem. We used the same network architecture in all experiments. The transition rates and growth rates were parameterized by two separate five layer neural networks, each with 128 hidden units per layer. The TV and RME were evaluated between the direct solution and the numerical solution and averaged over all evaluation time points. Since total variation does not account for differences in the distances between states, we additionally computed the normalized Wasserstein distance between the normalized direct solution and numerical solution at each evaluation time point, normalized by the maximum value of the edge weight \(C_{ij}\).
\[
\mathcal{W}(\widehat\mu,\mu)
=\frac{1}{\max_{ij}C_{ij}}
\min_{\gamma\in\Pi(\widehat\mu/|\widehat\mu|,\,\mu/|\mu|)}
\sum_{ij}C_{ij}\gamma_{ij},
\]
The average normalized Wasserstein distance over all evaluation time points were reported. 

\begin{table}[h]
\centering
\caption{Mean TV, RME, and normalized Wasserstein distance between closed form solution and CFM solution.}
\label{tab:cfm-direct-validation}
\begin{tabular}{lrrr}
\toprule
Dataset & Mean normalized TV & Mean RME & Mean $\mathcal{W}$ \\
\midrule
Synthetic (50 states)    & $3.79\times10^{-2}$ & $1.29\times10^{-3}$ & $1.07\times10^{-3}$ \\
Synthetic (1,000 states) & $1.14\times10^{-1}$ & $8.41\times10^{-3}$ & $9.89\times10^{-4}$ \\
EMCC (7 states)          & $8.30\times10^{-3}$ & $5.51\times10^{-3}$ & $1.41\times10^{-3}$ \\
\bottomrule
\end{tabular}
\end{table}

As shown in Table \ref{tab:cfm-direct-validation}, CFM reproduced the Direct measure flow with small RME in all three settings. The TV was below $0.04$ on the 50-state toy and EMCC and remained $0.114$ in the substantially larger 1,000-state problem.  Nevertheless, $\mathcal{W}$ remained close to $10^{-3}$ in all three datasets, showing that the residual distributional discrepancy was concentrated between nearby states. These results verify that CFM can recover the marginal WFR solution without explicitly constructing marginal training targets.

\subsection{Application to real biological data}
\label{sec:real-data-application}
We evaluated WFR-DFM on three single-cell datasets covering dynamic reconstruction, temporal interpolation, and robustness to annotation resolution. In all analyses, annotated cell types defined the discrete states, and distances between their mean gene expression profiles defined the WFR ground cost. Preprocessing and analysis settings are described in Supplementary Material~\ref{sec:real-data-analysis}.

\subsubsection{WFR-DFM reconstructs the dynamics of Mouse hematopoietic differentiation}
\label{sec:weinreb-application}
\begin{figure}[!t]
    \centering
    \includegraphics[
        width=0.9\textwidth,
        height=0.78\textheight,
        keepaspectratio
    ]{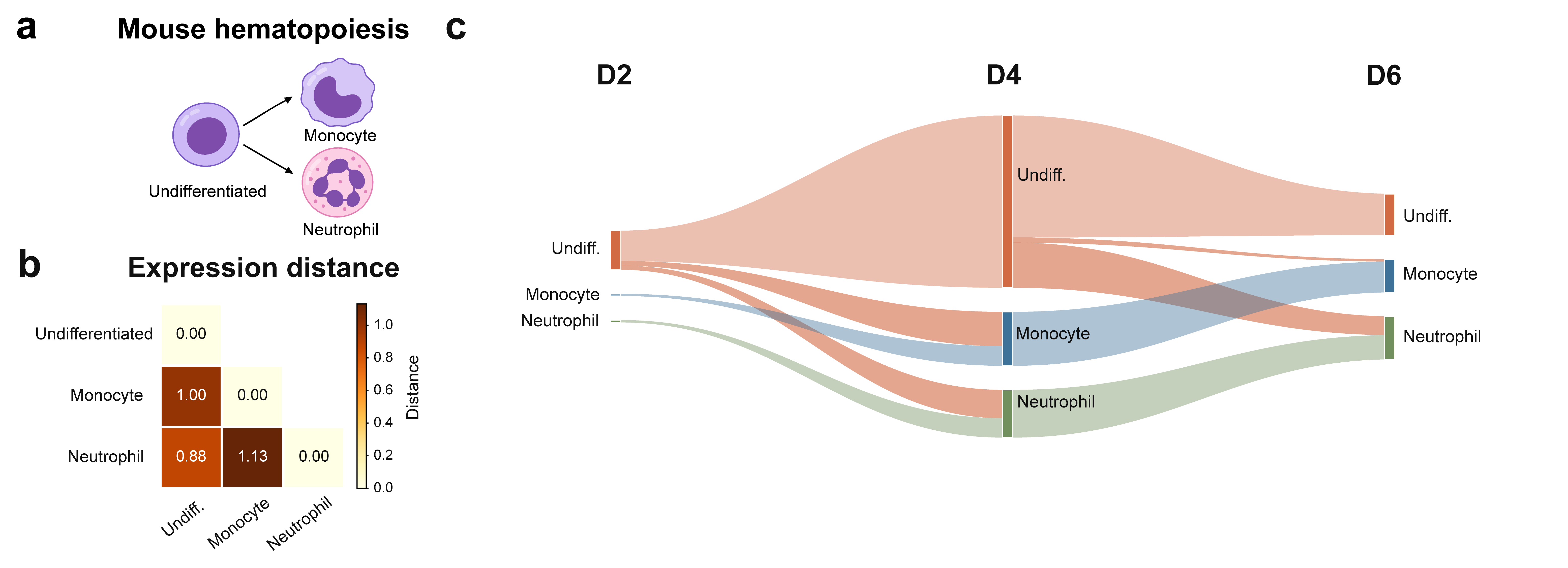}
    \caption{
        \textbf{Reconstructed dynamics of mouse hematopoietic differentiation.}
        \textbf{a}, Literature supported differentiation from undifferentiated
        progenitors toward monocytes and neutrophils \cite{weinreb2020lineage}.
        \textbf{b}, Pairwise expression distances between the three cell states used
        to define the WFR ground cost.
        \textbf{c}, Reconstructed WFR dynamics between consecutive observations.
        Node labels show sampling adjusted cell numbers. Ribbons represent inferred
        transitions, and changes in ribbon width reflect net population gain or loss.
    }
    \label{fig:Weinreb-wfr-dfm}
\end{figure}

We analysed a lineage traced mouse hematopoietic differentiation dataset
containing 96,371 cells sampled on days 2, 4 and 6
\cite{weinreb2020lineage}. We focused on the neutrophil and monocyte
differentiation branch and retained three annotated cell states:
undifferentiated cells, monocytes and neutrophils. Further dataset details are provided in
Supplementary Material~\ref{sec:app_hematopoiesis}. Figure~\ref{fig:Weinreb-wfr-dfm}a summarizes the differentiation from
undifferentiated progenitors toward monocytes and neutrophils supported by
lineage tracing \cite{weinreb2020lineage}, providing a biological reference
for comparison with the reconstructed dynamics.

The expression distance matrix shows that the distance between undifferentiated cells and either of the two mature cell states is smaller than the distance between monocytes and neutrophils (Figure \ref{fig:Weinreb-wfr-dfm}b). This is consistent with the biological understanding that undifferentiated cells can differentiate into Monocyte and Neutrophil, whereas Monocyte and Neutrophil cannot naturally transition into one another. WFR-DFM reconstructed the dynamics over day 2 to day 4 and day 4 to day 6 intervals (Fig.~\ref{fig:Weinreb-wfr-dfm}c). The inferred dynamics show that the undifferentiated cell population initially increases and continuously differentiates into monocytes and neutrophils over both time intervals, consistent with the two differentiation branches supported by existing biological knowledge\cite{weinreb2020lineage}. These results demonstrate that WFR-DFM can reconstruct biologically consistent continuous-time dynamics from sparsely observed snapshots.

\subsubsection{WFR-DFM interpolates the unobserved time point on developing mouse brain}

To further evaluate WFR-DFM quantitatively, we test WFR-DFM on a developing mouse brain dataset. It contains 120,576 cells collected at nine embryonic time points from E8.0 to E16.0 \cite{lamanno2021molecular}. Seven annotated neural cell types were retained, including neural tube cells, radial glia, neuroblasts, neurons, glioblasts, oligodendrocytes and ependymal cells. Dataset construction and interpolation settings are described in Supplementary Material~\ref{sec:app_dmb}. We evaluated temporal interpolation by holding out one time point at a time. For each excluded time point, only neighboring observations were used to construct the WFR bridge, and the inferred distribution was evaluated at the held-out time point.

\begin{figure}[!t]
    \centering
    \includegraphics[
        width=0.8\textwidth,
        keepaspectratio
    ]{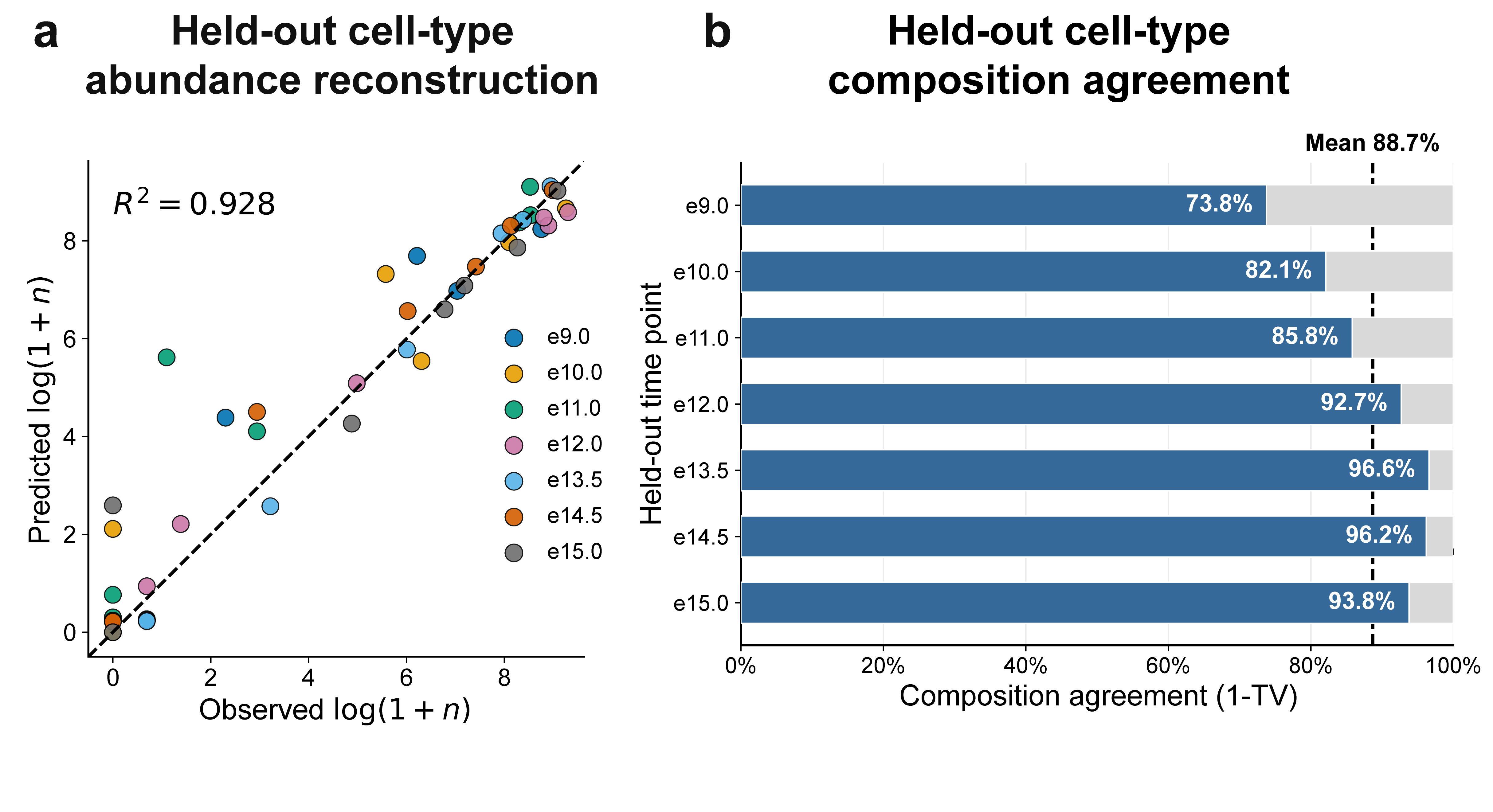}
    \caption{
    \textbf{Held-out interpolation of developmental dynamics in developing mouse brain.}
    \textbf{a}, Recovery of cell type abundances at excluded time points. Each color corresponds to one held-out time point, containing seven points representing the seven cell types. The \(x-\)axis shows the observed abundance and the \(y-\)axis shows the WFR predicted abundance. Values are shown as $\log(1+n)$, and the dashed line indicates exact agreement.
    \textbf{b}, Recovery of cell type composition at excluded time points. Each bar represents the composition agreement for one interpolation task, defined as $1-\mathrm{TV}$ between predicted and observed normalized compositions. The dashed line indicates the mean composition agreement across all interpolation tasks.
    }
    \label{fig:mouse-brain-wfr-dfm}
\end{figure}

We assessed interpolation accuracy from two complementary perspectives: recovery of cell type abundance and recovery of cell type composition. Fig.~\ref{fig:mouse-brain-wfr-dfm}a compares predicted and observed abundances for each cell type at held-out time points. The strong agreement on the logarithmic scale ($R^2=0.928$) indicates accurate recovery of cell type abundance at held-out time points.

In addition to abundance reconstruction, we evaluated whether the relative composition of cell types was preserved at held-out time points. Fig.~\ref{fig:mouse-brain-wfr-dfm}b summarizes the composition agreement, defined as $1-\mathrm{TV}$, across the seven interpolation tasks. The mean TV across these tasks is 0.113, corresponding to a mean composition agreement of 88.7\%. Together, these results show that WFR-DFM can interpolate developmental dynamics at unobserved time points by recovering both cell type abundance and composition from neighbouring observations.

\subsubsection{Consistency across annotation resolutions}
\label{sec:consistency}
\begin{figure}[!t]
    \centering
    \includegraphics[width=\textwidth]{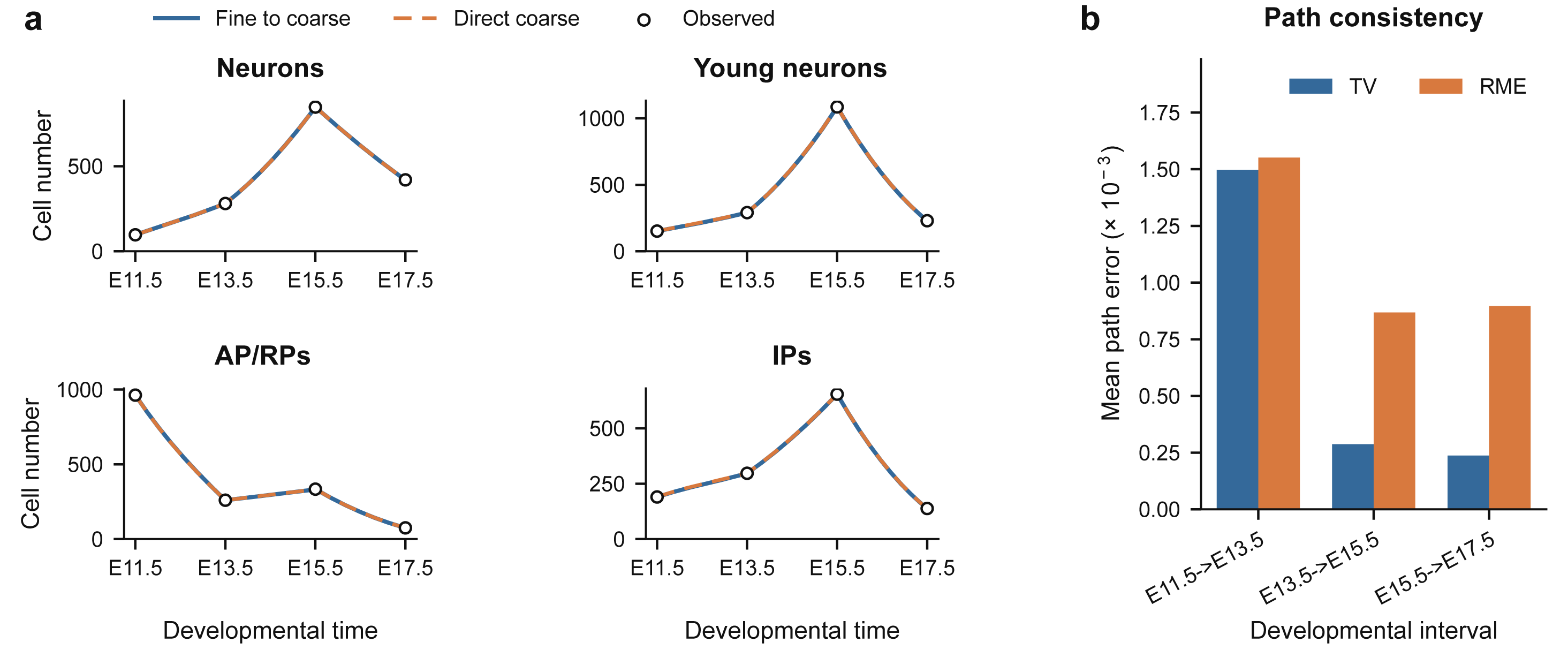}
    \caption{
    \textbf{Consistency of WFR trajectories across annotation resolutions.}
    \textbf{a}, Comparison of coarse trajectories obtained from two strategies: aggregating the trajectories inferred from seven fine cell types into four coarse cell classes (fine to coarse), and directly fitting WFR using the four coarse cell classes (direct coarse). Each panel shows one coarse cell class. Open circles indicate the observed coarse abundances, whereas solid and dashed curves represent the two reconstructed trajectories.
    \textbf{b}, Mean path errors for each developmental interval between the two reconstruction strategies. Errors are evaluated over the intermediate positions of each WFR bridge and quantified using total variation (TV) and relative mass error (RME).
    }
    \label{fig:emcc-coarsening-consistency}
\end{figure}

Cell annotations are often available at different resolutions, ranging from fine subtypes to broader cell classes. We evaluated the sensitivity of WFR inference to annotation resolution using the embryonic mouse cerebral cortex (EMCC) dataset \cite{EMCC}, which contains 6,316 cells collected at four time points (E11.5, E13.5, E15.5 and E17.5) and annotated at two resolutions. At the fine resolution, the data is annotated into seven cell types spanning progenitor, intermediate progenitor and neuronal states. At coarse resolution, it is annotated into four. We compared two analysis strategies. In the first, WFR trajectories were inferred using the seven annotated cell types and then aggregated into four coarse cell classes. In the second, the observations were first aggregated into the same four classes, and WFR was fitted directly at the coarse resolution. Dataset preprocessing and the detailed aggregation procedure are provided in Supplementary Material~\ref{sec:app_emcc}.
The two strategies produced highly consistent coarse trajectories across developmental time (Fig.~\ref{fig:emcc-coarsening-consistency}a). The fine to coarse trajectories closely overlapped with those obtained from direct coarse inference, while both approaches matched the observed coarse abundances. 

To quantify the agreement between the two inferred trajectories, we evaluated the discrepancy at 201 equally spaced path positions for each WFR bridge. The two shared endpoints were excluded, and the mean path errors were computed over the remaining intermediate positions:
\begin{equation*}
    \mathrm{Mean\ TV}
    =
    \frac{1}{N}
    \sum_{k=1}^{N}
    \mathrm{TV}
    \left(
        \frac{\mu_k^{\mathrm{fine\rightarrow coarse}}}{\left|\mu_k^{\mathrm{fine\rightarrow coarse}}\right|},
        \frac{\mu_k^{\mathrm{direct\ coarse}}}{\left|\mu_k^{\mathrm{direct\ coarse}}\right|}
    \right),
\end{equation*}
where $N=199$ denotes the number of intermediate path positions. Here, $\mu_k^{\mathrm{fine\rightarrow coarse}}$ and $\mu_k^{\mathrm{direct\ coarse}}$ denote the coarse state measures at the $k$th intermediate path position obtained from the two strategies, respectively. The mean RME was computed analogously. The resulting discrepancies remained small across all developmental intervals (Fig.~\ref{fig:emcc-coarsening-consistency}b), with mean path errors on the order of $10^{-3}$. These results suggest that, under this aggregation scheme, WFR inference is robust to annotation resolution.

\section{Conclusion}\label{conclusion}
This paper proposes a dynamic discrete Wasserstein--Fisher--Rao (WFR) framework for reconstructing unbalanced population dynamics on discrete state spaces. Building on the travelling Dirac technique, we construct a dynamic realization of the static WFR problem. We further derive an explicit unbalanced Markov representation of WFR measure flows, in which state-dependent transition rates and growth rates jointly characterize state transitions and mass variation. When direct marginalization is intractable, unbalanced discrete conditional flow matching provides a scalable approximation by learning tractable conditional targets. We also characterize the behavior of the framework under state aggregation, establishing exact consistency for block-homogeneous costs and stability under small geometric perturbations. Numeric experiments on synthetic data show that explicitly modeling growth prevents changes in mass from being spuriously interpreted as state transitions. Applications to single-cell transcriptomic data further demonstrate biologically consistent dynamical reconstruction, temporal interpolation, and robustness across different annotation resolutions.

Future work may extend the multiscale robustness analysis to continuous spaces, with the goal of establishing a connection between the continuous limit of WFR flows on discretized state spaces and WFR flows in continuous spaces. Such a theory could provide a foundation for developing efficient multiscale algorithms for WFR dynamics in continuous state spaces.


\clearpage
\appendix
\section{Dataset details}
\label{app:data}
\subsection{Synthetic datasets}
\subsubsection{50-state synthetic data generation}
The 50-state synthetic data consists of states
$\{x_i\}_{i=1}^{50}$ sampled independently from $\mathcal{U}[0,1]^2$. The edge weights are calculated as $C_{ij}=\|x_i-x_j\|_2^2$. At the four observation time points, the observed measures are supported on
\[
S_0=\{1,\ldots,18\},\quad
S_1=\{9,\ldots,30\},\quad
S_2=\{21,\ldots,43\},\quad
S_3=\{33,\ldots,50\},
\]
so consecutive distributions have partially overlapping supports. For
$i\in S_k$, the mass is defined as
\[
\mu_k(i)
=
\ell_k+
\exp\left[-\frac{1}{2}
\left(\frac{i-\bar i_k}{5.5}\right)^2\right]
+\varepsilon_{k,i},
\qquad
(\ell_0,\ell_1,\ell_2,\ell_3)=(0.5,0.8,1.1,1.4),
\]
where $\bar i_k$ is the midpoint of $S_k$ and
$\varepsilon_{k,i}\overset{i.i.d}{\sim}\mathcal{U}[0,0.15]$. Thus, the synthetic data exhibits both support displacement and changes in total mass.

\subsubsection{1000-state synthetic data generation}
The 1000-state space is obtained from a regular $32\times32$ grid on $[0,1]^2$. Specifically, the grid points are ordered row by row as
\[
x_{32j+i+1}=\left(\frac{i}{31},\frac{j}{31}\right),
\qquad i,j\in\{0,\ldots,31\},
\]
where the horizontal index $i$ varies first. We retain the first 1000 points
under this ordering. Equivalently, the resulting state space is
\[
\mathcal X
=
\left\{
\left(\frac{i}{31},\frac{j}{31}\right):
i=0,\ldots,31,\ j=0,\ldots,30
\right\}
\cup
\left\{
\left(\frac{i}{31},1\right):
i=0,\ldots,7
\right\}.
\]
Thus, it contains all 992 points in the first 31 grid rows and the first eight
points in the final row. The edge weights are also calculated as $C_{ij}=\|x_i-x_j\|_2^2$. At observation time $k$, its unnormalized density is defined as
\[
\rho_k(x)
=
0.025+
\exp\left(-\frac{\|x-c_k\|_2^2}{2(0.115)^2}\right)
+
0.38\exp\left(-\frac{\|x-\widetilde c_k\|_2^2}{2(0.18)^2}\right),
\]
where
\[
\begin{aligned}
(c_0,c_1,c_2,c_3)
&=((0.14,0.18),(0.38,0.72),(0.67,0.33),(0.86,0.76)),\\
(\widetilde c_0,\widetilde c_1,\widetilde c_2,\widetilde c_3)
&=((0.78,0.84),(0.69,0.24),(0.30,0.77),(0.18,0.28)).
\end{aligned}
\]
The discrete measure is defined as
\[
\mu_k(i)=M_k\frac{\rho_k(x_i)}
{\sum_{j=1}^{1000}\rho_k(x_j)},
\qquad
(M_0,M_1,M_2,M_3)=(1.00,1.35,0.72,1.18).
\]
The positive background gives full support, while the moving Gaussian components and prescribed total masses also induce simultaneous transport and
growth or decay.


\subsection{Real-data analysis}
\label{sec:real-data-analysis}

\subsubsection{Mouse hematopoietic differentiation with lineage tracing}
\label{sec:app_hematopoiesis}

The in vitro mouse hematopoietic differentiation dataset was reported by
Weinreb et al. \cite{weinreb2020lineage}. We analysed the published
neutrophil and monocyte differentiation trajectory, containing 96,371 cells
sampled on days 2, 4 and 6. Cells were grouped into three annotated states:
undifferentiated cells, monocytes and neutrophils. We used the provided
total count normalized expression matrix and applied a $\log(1+x)$
transformation. A total of 3,000 highly variable genes were then selected
using the Seurat style dispersion criterion.

Following the sampling correction of Sha et al. \cite{TIGON}, we constructed
the cell-state masses as
\begin{equation}
    \rho_t(i)=\frac{n_t(i)}{s_t\sum_j n_2(j)},
    \qquad
    s_2=1,\quad s_4=0.30,\quad s_6=0.70,
\end{equation}
where $n_t(i)$ is the observed number of cells in state $i$ at day $t$,
and $s_t$ denotes the sampling fraction at that time point.
The resulting quantities define sampling adjusted masses relative to the
total sequenced cell count at day 2.

All three time points were used for HVG selection and ground cost construction.
WFR bridges were then computed for the day 2 to day 4 and day 4 to day 6
intervals using the sampling adjusted masses. We used $\delta=1.0$ for the
reported analysis.

\subsubsection{Developing mouse brain}
\label{sec:app_dmb}

We used the developing mouse brain atlas of La Manno et al. \cite{lamanno2021molecular}. The analysis included a subset of 120,576 cells sampled at nine embryonic time points: E8.0, E9.0, E10.0, E11.0, E12.0, E13.5, E14.5, E15.0 and E16.0. Expression counts were normalized to $10^4$ counts per cell and transformed using $\log(1+x)$. For this dataset, 3,000 HVGs were selected using the Seurat style dispersion criterion. Seven annotated neural cell types were retained: neural tube cells, radial glia, neuroblasts, neurons, glioblasts, oligodendrocytes and ependymal cells.

Temporal interpolation was evaluated by excluding each of the seven internal time points in turn. In each interpolation fold, cells from the excluded time point were removed before HVG selection, centroid estimation and ground cost construction. A WFR bridge was computed between the immediately preceding and following observed time points and evaluated at

\begin{equation*}
    \tau_h
    =
    \frac{t_h-t_{h-1}}
    {t_{h+1}-t_{h-1}},
\end{equation*}

where $t_h$ denotes the excluded time point. HVG selection, centroid estimation and distance normalization were performed independently within each fold using only the retained observations. We used $\delta=0.70$ for all developing mouse brain analyses.

\subsubsection{Embryonic mouse cerebral cortex}
\label{sec:app_emcc}

The embryonic mouse cerebral cortex dataset was reported by Yuzwa et al. \cite{EMCC}. We used the processed version distributed with Tempora \cite{tran2020tempora}, comprising 6,316 cells collected at four time points: E11.5, E13.5, E15.5 and E17.5. The provided normalized expression matrix was used, and negligible negative numerical values were set to zero. For this dataset, 3,000 highly variable genes (HVGs) were selected using the Seurat style dispersion criterion. The cells were annotated into seven cell types, including one neuron type, two young neuron types, two apical progenitor or radial precursor (AP/RP) types, and two intermediate progenitor (IP) types.

All four time points were included in HVG selection, cost construction and computation of the three successive WFR bridges. We used $\delta=1.0$ for the reported analysis.

For the annotation resolution analysis, the seven fine cell types were
aggregated into four coarse cell classes: neurons, young neurons, AP/RPs
and IPs. For the fine to coarse strategy, the mass of each coarse class at
every path position was obtained by summing the masses of its constituent
fine cell types. For direct coarse inference, coarse class centroids were
computed as cell number weighted averages of the corresponding fine type
centroids, and the WFR problem was solved directly on the four coarse
classes. TV and RME between the two trajectories were evaluated at the 199
intermediate positions of each 201 point path, excluding the shared observed
endpoints.

\subsubsection{WFR construction}
\label{sec:PandWFR}

For all datasets, the provided cell type annotations defined the discrete
states. The WFR ground cost was constructed from gene expression differences
between cell types. Let $\bar{\boldsymbol z}_i$ denote the mean HVG expression
vector of cell type $i$, calculated using all cells retained in the current
analysis. The pairwise expression distance was defined as

\begin{equation*}
    d_{ij}^{\mathrm{raw}}
    =
    \left\|
        \bar{\boldsymbol z}_i-\bar{\boldsymbol z}_j
    \right\|_2,
    \qquad
    s_d
    =
    \operatorname{median}_{i<j,\,d_{ij}^{\mathrm{raw}}>0}
    d_{ij}^{\mathrm{raw}},
    \qquad
    C_{ij}
    =
    \left(
        \frac{d_{ij}^{\mathrm{raw}}}{s_d}
    \right)^2.
\end{equation*}

At each observation time, cell type abundances were used to construct the
corresponding discrete measure. For each WFR bridge, the static WFR
semi-coupling was computed between the two endpoint observations using the
ground cost $C$. The resulting bridges provided time dependent cell type
masses, transition rates and growth rates.



\section{\textbf{Proofs}}
\label{app:deferred-proofs}

\subsection{\textbf{Proof of Proposition~\ref{prop:ot metric}}}
\label{app:proof-ot-metric}

\begin{proof}
To prove the metric property, it is sufficient to show that \(\text{OT}^{\frac{1}{p}}\) satisfies the triangle inequality. 
\begin{equation}
\text{OT}^{\frac{1}{p}}(\alpha,\beta)+\text{OT}^{\frac{1}{p}}(\beta,\gamma)\ge\text{OT}^{\frac{1}{p}}(\alpha,\gamma)
\end{equation}
where \(\alpha,\beta,\gamma\) are three distributions over \(\mathcal{X}\).  

Let \(\pi_{\alpha\beta},\pi_{\beta\gamma}\) be the optimal transport coupling between \(\alpha,\beta\), and \(\beta,\gamma\). A gluing is defined as 
\begin{equation}
(\pi_{\alpha\beta\gamma})_{ijk}=\left\{
\begin{aligned}
&(\pi_{\alpha\beta})_{ij}(\pi_{\beta\gamma})_{jk}/\beta_j,\quad &\beta_j\neq0\\
&0,\quad&\beta_j=0
\end{aligned}
\right .
\end{equation}

One can easily check that 
\[
\sum_{k}(\pi_{\alpha\beta\gamma})_{ijk}=(\pi_{\alpha\beta})_{ij},\quad\sum_{i}(\pi_{\alpha\beta\gamma})_{ijk}=(\pi_{\beta\gamma})_{jk}
\]

Based on the gluing, a coupling between \(\alpha\) and \(\gamma\) can be constructed as \((\pi_{\alpha\gamma})_{ik}=\sum_{j}(\pi_{\alpha\beta\gamma})_{ijk}\). It is a coupling between \(\alpha\) and \(\gamma\), since
\begin{equation*}
\begin{aligned}
\sum_{i}(\pi_{\alpha\gamma})_{ik}
&=\sum_{i,j}(\pi_{\alpha\beta\gamma})_{ijk}
=\sum_{j}(\pi_{\beta\gamma})_{jk}=\gamma_k,\\
\sum_{k}(\pi_{\alpha\gamma})_{ik}
&=\sum_{j,k}(\pi_{\alpha\beta\gamma})_{ijk}
=\sum_{j}(\pi_{\alpha\beta})_{ij}=\alpha_i.
\end{aligned}
\end{equation*}

Therefore, by the metric property of \(C^{\frac{1}{p}}\) and the Minkovski inequality, we have
\begin{equation}
\begin{aligned}
\text{OT}^{\frac{1}{p}}(\alpha,\gamma)
&\le\left(\sum_{i,k}C_{ik}(\pi_{\alpha\gamma})_{ik}\right)^{\frac{1}{p}}\\
&=\left(\sum_{i,j,k}C_{ik}(\pi_{\alpha\beta\gamma})_{ijk}\right)^{\frac{1}{p}}\\
&\le\left(\sum_{i,j,k}(C_{ij}^{\frac{1}{p}}+C_{jk}^{\frac{1}{p}})^p
(\pi_{\alpha\beta\gamma})_{ijk}\right)^{\frac{1}{p}}\\
&\le\left(\sum_{i,j,k}C_{ij}(\pi_{\alpha\beta\gamma})_{ijk}\right)^{\frac{1}{p}}
+\left(\sum_{i,j,k}C_{jk}(\pi_{\alpha\beta\gamma})_{ijk}\right)^{\frac{1}{p}}\\
&=\left(\sum_{i,j}C_{ij}(\pi_{\alpha\beta})_{ij}\right)^{\frac{1}{p}}
+\left(\sum_{j,k}C_{jk}(\pi_{\beta\gamma})_{jk}\right)^{\frac{1}{p}}\\
&=\text{OT}^{\frac{1}{p}}(\alpha,\beta)
+\text{OT}^{\frac{1}{p}}(\beta,\gamma).
\end{aligned}
\end{equation}

When \(p=1\), let \(\gamma(t)=(1-t)\alpha+t\beta\). Consider the transport plan \(\pi_{st}=(t-s)\pi_{\alpha\beta}+(1-t)\text{diag}(\alpha)+s\cdot\text{diag}(\beta)\) which transports \(\gamma(s)\) to \(\gamma(t)\). By direct calculation
\begin{equation}
\text{OT}(\gamma(s),\gamma(t))\le\sum_{i,j}C_{ij}(\pi_{st})_{ij}=(t-s)\text{OT}(\alpha,\beta)\quad(s<t)
\end{equation}

By the triangle inequality which we just proved, 
\begin{equation}
\begin{aligned}
\text{OT}(\alpha,\beta)
&\le\text{OT}(\alpha,\gamma(s))+\text{OT}(\gamma(s),\gamma(t))
+\text{OT}(\gamma(t),\beta)\\
&\le [s+(t-s)+(1-t)]\text{OT}(\alpha,\beta)\\
&=\text{OT}(\alpha,\beta).
\end{aligned}
\end{equation}
Therefore, all the equalities must hold. We have \(\text{OT}(\gamma(s),\gamma(t))=(t-s)\text{OT}(\alpha,\beta)\), i.e. \(\gamma(t)=(1-t)\alpha+t\beta\) is the geodesic under the OT metric.
\end{proof}

\subsection{\textbf{Proof of Proposition~\ref{prop:existence}}}
\label{app:proof-existence}

\begin{proof}
The feasible set
\begin{equation*}
    \left\{(a,b)\in\Rnum_{\ge0}^{N\times N}\times\Rnum_{\ge0}^{N\times N}:
    \sum_j a_{ij}=\rho_i,\ \sum_i b_{ij}=\eta_j\right\}
\end{equation*}
is closed and bounded in a finite-dimensional Euclidean space, hence compact. The objective in \eqref{eq:static-wfr} is continuous on this compact set, and therefore attains its minimum. If $\eta=0$, the constraints force $b=0$ and the objective becomes $4\delta^2\sum_{i,j}a_{ij}=4\delta^2|\rho|$. The case $\rho=0$ is identical.
\end{proof}

\subsection{\textbf{Proof of Proposition~\ref{prop:basic-properties}}}
\label{app:proof-basic-properties}

\begin{proof}
Nonnegativity follows from the pointwise identity
\begin{equation*}
    a+b-2\zeta\sqrt{ab}
    =
    (\sqrt a-\sqrt b)^2+2(1-\zeta)\sqrt{ab}\geq0,
    \qquad 0\leq \zeta\leq1.
\end{equation*}
Symmetry follows by exchanging the indices $i,j$ and the variables $a,b$. Homogeneity follows by replacing $(a,b)$ with $(\lambda a,\lambda b)$.

For the upper bound in \eqref{eq:mass-bounds}, dropping the negative terms
\[
-2\sqrt{a_{ij}b_{ij}}\,
\cos\left(\frac{\sqrt{C_{ij}}}{2\delta}\right)
\]
yields:
\begin{equation*}
    \WFR_\delta^2(\rho,\eta)
    \leq 4\delta^2(|\rho|+|\eta|).
\end{equation*}
This implies the displayed upper bound, which is weaker but simpler. For the lower bound, write $\zeta_{ij}=\cos(\sqrt{C_{ij}}/2\delta)\leq1$. For any feasible $(a,b)$,
\begin{equation*}
\begin{split}
    \sum_{i,j}\left(a_{ij}+b_{ij}-2\zeta_{ij}\sqrt{a_{ij}b_{ij}}\right)
    &\geq
    \sum_{i,j}\left(a_{ij}+b_{ij}-2\sqrt{a_{ij}b_{ij}}\right)\\
    &\geq
    |\rho|+|\eta|-2\sqrt{|\rho||\eta|}
    =
    \left(\sqrt{|\rho|}-\sqrt{|\eta|}\right)^2,
\end{split}
\end{equation*}
where the second inequality is Cauchy's inequality. Multiplying by $4\delta^2$ yields the lower bound.

If $\rho=\eta$, choosing $a_{ii}=b_{ii}=\rho_i$ and all off-diagonal entries zero gives zero cost. Conversely, if the cost is zero, every active pair must satisfy $a_{ij}=b_{ij}$ and $\zeta_{ij}=1$ whenever $a_{ij}b_{ij}>0$. Since $\zeta_{ij}=1$ implies $C_{ij}=0$, hence $i=j$. Therefore \(a_{ij}=b_{ij}=0, i\neq j;a_{ii}=b_{ii}\), which yields \(\rho=\eta\).
\end{proof}

\subsection{\textbf{Proof of Theorem~\ref{thm:metric}}}
\label{app:proof-wfr-metric}

\begin{proof}
It is sufficient to prove the triangle inequality. Since \(1\le p\le 2\), the assumption that \(C_{ij}^{1/p}\) is a metric implies that \(C_{ij}^{1/2}\) is also a metric. Indeed, for any \(i,j,k\),
\[
    C_{ik}^{1/p}\le C_{ij}^{1/p}+C_{jk}^{1/p}.
\]
Raising both sides to the power \(p/2\le 1\), we have
\[
    C_{ik}^{1/2}
    \le
    \left(C_{ij}^{1/p}+C_{jk}^{1/p}\right)^{p/2}
    \le
    C_{ij}^{1/2}+C_{jk}^{1/2}.
\]
Hence \(C^{1/2}\) satisfies the triangle inequality.

Since \(\max_{i,j}C_{ij}<\pi^2\delta^2\), we use the static semi-coupling formulation
\begin{equation}
\begin{aligned}
    \WFR_\delta(\rho,\eta)
    &=
    2\delta
    \inf_{(A_0,A_1)\in\Gamma(\rho,\eta)}
    \Biggl[
    \sum_{i,j}
    (A_0)_{ij}+(A_1)_{ij}
    \\
    &\qquad
    -2\sqrt{(A_0)_{ij}(A_1)_{ij}}
    \cos\left(\frac{\sqrt{C_{ij}}}{2\delta}\right)
    \Biggr]^{1/2},
\end{aligned}
\end{equation}
Let \(\alpha,\beta,\gamma\in\mathcal M_+(\mathcal X)\). Take optimal semi-couplings \((A_0,A_1)\in\Gamma(\alpha,\beta),(B_0,B_1)\in\Gamma(\beta,\gamma).\) Since \(\sum_i (A_1)_{ij}=\beta_j=\sum_k (B_0)_{jk},\) we glue \(A_1\) and \(B_0\) along the common marginal \(\beta\). Define
\[
    \Pi_{ijk}
    =
    \begin{cases}
        \dfrac{(A_1)_{ij}(B_0)_{jk}}{\beta_j},
        & \beta_j\neq 0,\\[1.2ex]
        0,
        & \beta_j=0.
    \end{cases}
\]
Then \(\sum_k \Pi_{ijk}=(A_1)_{ij},\sum_i \Pi_{ijk}=(B_0)_{jk}.\) Lift \(A_0\) and \(B_1\) to the same gluing by setting
\[
    (\Pi_0)_{ijk}
    =
    \begin{cases}
        \dfrac{(A_0)_{ij}(B_0)_{jk}}{\beta_j},
        & \beta_j\neq 0,\\[1.2ex]\dfrac{(A_0)_{ij}}{N},
        & \beta_j=0,
    \end{cases}
    \qquad
    (\Pi_1)_{ijk}
    =
    \begin{cases}
        \dfrac{(A_1)_{ij}(B_1)_{jk}}{\beta_j},
        & \beta_j\neq 0,\\[1.2ex]\dfrac{(B_1)_{jk}}{N},
        & \beta_j=0.
    \end{cases}
\]
One can easily check that
\[
\sum_{j,k}(\Pi_0)_{ijk}=\alpha_i,
\qquad
\sum_{i,j}(\Pi_1)_{ijk}=\gamma_k.
\]
Therefore,
\[
\left(\sum_j(\Pi_0)_{ijk},\sum_j(\Pi_1)_{ijk}\right)_{i,k}
\]
defines a semi-coupling between \(\alpha\) and \(\gamma\). For simplicity, denote
\[
d_{ij}:=\frac{\sqrt{C_{ij}}}{2\delta}.
\]
For each triple \((i,j,k)\), define
\[
\begin{aligned}
    X_{ijk}
    &=
    \left[
    (\Pi_0)_{ijk}+\Pi_{ijk}
    -2\sqrt{(\Pi_0)_{ijk}\Pi_{ijk}}\cos d_{ij}
    \right]^{1/2},\\
    Y_{ijk}
    &=
    \left[
    \Pi_{ijk}+(\Pi_1)_{ijk}
    -2\sqrt{\Pi_{ijk}(\Pi_1)_{ijk}}\cos d_{jk}
    \right]^{1/2}.
\end{aligned}
\]
By the Minkowski inequality,
\begin{equation}
\begin{aligned}
    \frac{1}{2\delta}
    \left[
    \WFR_\delta(\alpha,\beta)+\WFR_\delta(\beta,\gamma)
    \right]
    &=
    \left(\sum_{i,j,k}X_{ijk}^2\right)^{1/2}
    +
    \left(\sum_{i,j,k}Y_{ijk}^2\right)^{1/2}\\
    &\ge
    \left(\sum_{i,j,k}(X_{ijk}+Y_{ijk})^2\right)^{1/2}.
\end{aligned}
\end{equation}

Now we use the inequality
\begin{equation}
\label{eq:cone-triangle}
\left[
a+b-2\sqrt{ab}\cos d_1
\right]^{1/2}
+
\left[
b+c-2\sqrt{bc}\cos d_2
\right]^{1/2}\ge
\left[
a+c-2\sqrt{ac}\cos(d_1+d_2)
\right]^{1/2},
\end{equation}
valid for \(a,b,c\ge0\) and \(d_1,d_2\in[0,\pi/2)\). This is simply the triangle inequality in the Euclidean plane. Applying \eqref{eq:cone-triangle} with
\(a=(\Pi_0)_{ijk},b=\Pi_{ijk},c=(\Pi_1)_{ijk},\) we get
\[
    X_{ijk}+Y_{ijk}
    \ge
    \left[
    (\Pi_0)_{ijk}+(\Pi_1)_{ijk}
    -2\sqrt{(\Pi_0)_{ijk}(\Pi_1)_{ijk}}
    \cos(d_{ij}+d_{jk})
    \right]^{1/2}.
\]
Since \(C^{1/2}\) is a metric,
\[
    d_{ik}
    =
    \frac{\sqrt{C_{ik}}}{2\delta}
    \le
    \frac{\sqrt{C_{ij}}+\sqrt{C_{jk}}}{2\delta}
    =
    d_{ij}+d_{jk}.
\]
which further gives \(\cos(d_{ik})\ge\cos(d_{ij}+d_{jk}).\) Hence
\[
\begin{aligned}
    X_{ijk}+Y_{ijk}
    &\ge
    \left[
    (\Pi_0)_{ijk}+(\Pi_1)_{ijk}
    -2\sqrt{(\Pi_0)_{ijk}(\Pi_1)_{ijk}}
    \cos(d_{ik})
    \right]^{1/2}.
\end{aligned}
\]
Therefore,
\begin{equation}
\begin{aligned}
    \frac{1}{2\delta}
    \left[
    \WFR_\delta(\alpha,\beta)+\WFR_\delta(\beta,\gamma)
    \right]
    &\ge
    \Biggl[
    \sum_{i,j,k}\Bigl(
    (\Pi_0)_{ijk}+(\Pi_1)_{ijk}
    \\
    &\qquad
    -2\sqrt{(\Pi_0)_{ijk}(\Pi_1)_{ijk}}
    \cos(d_{ik})\Bigr)
    \Biggr]^{1/2}.
\end{aligned}
\end{equation}

Finally, by Cauchy's inequality,
\[
    \sum_j
    \sqrt{(\Pi_0)_{ijk}(\Pi_1)_{ijk}}
    \le
    \sqrt{
    \left(\sum_j(\Pi_0)_{ijk}\right)
    \left(\sum_j(\Pi_1)_{ijk}\right)
    }.
\]
Thus
\begin{equation}
\begin{aligned}
&
\sum_{i,j,k}\Big(
(\Pi_0)_{ijk}+(\Pi_1)_{ijk}
-2\sqrt{(\Pi_0)_{ijk}(\Pi_1)_{ijk}}
\cos(d_{ik})\Big)\\
&\qquad\ge
\sum_{i,k}\Big(
\sum_j(\Pi_0)_{ijk}
+
\sum_j(\Pi_1)_{ijk}
-2
\sqrt{
\left(\sum_j(\Pi_0)_{ijk}\right)
\left(\sum_j(\Pi_1)_{ijk}\right)
}
\cos(d_{ik})\Big).
\end{aligned}
\end{equation}
Since the pair
\[
    \left(
    \sum_j(\Pi_0)_{ijk},
    \sum_j(\Pi_1)_{ijk}
    \right)_{i,k}
\]
is a semi-coupling between \(\alpha\) and \(\gamma\), by the definition of \(\WFR_\delta(\alpha,\gamma)\), we have
\[
    \frac{1}{2\delta}
    \left[
    \WFR_\delta(\alpha,\beta)+\WFR_\delta(\beta,\gamma)
    \right]
    \ge
    \frac{1}{2\delta}\WFR_\delta(\alpha,\gamma).
\]
Therefore,
\[
    \WFR_\delta(\alpha,\beta)+\WFR_\delta(\beta,\gamma)
    \ge
    \WFR_\delta(\alpha,\gamma).
\]
This proves the triangle inequality, and hence \(\WFR_\delta\) is a metric on \(\mathcal M_+(\mathcal X)\).
\end{proof}

\subsection{\textbf{Proof of Proposition~\ref{prop:entropy}}}
\label{app:proof-entropy}

\begin{proof}
Denote \(\omega_{ij}:=\cos\left(\frac{\sqrt{C_{ij}}}{2\delta}\right).\) Since \(\max_{i,j}C_{ij}<\pi^2\delta^2\), we have \(0<\omega_{ij}\le 1.\) By the static semi-coupling formulation of WFR,
\begin{equation}\label{eq:static-wfr-proof}
\begin{aligned}
    \frac{1}{4\delta^2}\WFR_\delta^2(\rho,\eta)
    =
    \inf_{a,b\ge0}
    \sum_{i,j}
    \left[
        a_{ij}+b_{ij}
        -2\sqrt{a_{ij}b_{ij}}\,\omega_{ij}
    \right],
\end{aligned}
\end{equation}
where the infimum is taken under the constraints \(\sum_j a_{ij}=\rho_i,\sum_i b_{ij}=\eta_j.\) We first derive the dual problem of \eqref{eq:static-wfr-proof}. Introduce Lagrange multipliers
\(\phi_i,\psi_j\) for the marginal constraints. The Lagrangian is
\[
\begin{aligned}
    L(a,b,\phi,\psi)
    &=
    \sum_{i,j}
    \left[
        a_{ij}+b_{ij}
        -2\sqrt{a_{ij}b_{ij}}\,\omega_{ij}
    \right]\\
    &\quad+\sum_i\phi_i\left(\rho_i-\sum_j a_{ij}\right)\\
    &\quad+\sum_j\psi_j\left(\eta_j-\sum_i b_{ij}\right).
\end{aligned}
\]
For fixed \(\phi,\psi\), minimizing with respect to \(a_{ij},b_{ij}\ge0\) gives a finite value if and only if
\[
    \phi_i\le 1,\qquad \psi_j\le 1,
    \qquad
    (1-\phi_i)(1-\psi_j)\ge \omega_{ij}^2,
\]
since the quadratic form in \(\sqrt{a_{ij}}\) and \(\sqrt{b_{ij}}\) is nonnegative precisely under this condition. Therefore, the dual problem is
\begin{equation}\label{eq:wfr-dual-phi}
\begin{aligned}
    \sup_{\phi,\psi}
    \left\{
        \sum_i\phi_i\rho_i+\sum_j\psi_j\eta_j
    \right\},
    \qquad
    \text{s.t. }
    \phi_i\le1,\ \psi_j\le1,\ 
    (1-\phi_i)(1-\psi_j)\ge \omega_{ij}^2.
\end{aligned}
\end{equation}

Now, introduce the change of variables \(u_i=-\log(1-\phi_i),v_j=-\log(1-\psi_j).\) Then \(\phi_i=1-e^{-u_i},\psi_j=1-e^{-v_j},\) and the constraint becomes \(u_i+v_j\le -2\log\omega_{ij}=c_{ij}.\) Thus, \eqref{eq:wfr-dual-phi} is equivalent to
\begin{equation}\label{eq:wfr-dual-uv}
    \sup_{u,v}
    \left\{
        \sum_i \rho_i(1-e^{-u_i})
        +
        \sum_j \eta_j(1-e^{-v_j})
    \right\},
    \qquad
    \text{s.t. } u_i+v_j\le c_{ij}.
\end{equation}

It remains to identify the primal problem associated with \eqref{eq:wfr-dual-uv}. Introduce a nonnegative Lagrange multiplier \(\gamma_{ij}\ge0\) for the constraint \(u_i+v_j\le c_{ij}\). Then \eqref{eq:wfr-dual-uv} has dual
\[
    \inf_{\gamma\ge0}
    \sup_{u,v}
    \left\{
        \sum_i \rho_i(1-e^{-u_i})
        +
        \sum_j \eta_j(1-e^{-v_j})
        +
        \sum_{i,j}\gamma_{ij}(c_{ij}-u_i-v_j)
    \right\}.
\]
Let
\[
    r_i=(\gamma\one)_i=\sum_j\gamma_{ij},
    \qquad
    s_j=(\gamma^T\one)_j=\sum_i\gamma_{ij}.
\]
Then the supremum over \(u,v\) separates as
\[
\begin{aligned}
    \sup_{u,v}
    &\left\{
        \sum_i \rho_i(1-e^{-u_i})-r_i u_i
        +
        \sum_j \eta_j(1-e^{-v_j})-s_j v_j
        +
        \sum_{i,j}c_{ij}\gamma_{ij}
    \right\}\\
    &=
    \sum_{i,j}c_{ij}\gamma_{ij}
    +
    \sum_i
    \sup_{u_i}
    \left[
        \rho_i(1-e^{-u_i})-r_i u_i
    \right]
    +
    \sum_j
    \sup_{v_j}
    \left[
        \eta_j(1-e^{-v_j})-s_j v_j
    \right].
\end{aligned}
\]
For each \(i\), by calculating the derivatives of \(u_i\), the maximizer satisfies \(\rho_i e^{-u_i}=r_i,\) and hence
\[
    \sup_{u_i}
    \left[
        \rho_i(1-e^{-u_i})-r_i u_i
    \right]
    =
    r_i\log\left(\frac{r_i}{\rho_i}\right)-r_i+\rho_i.
\]
Similarly,
\[
    \sup_{v_j}
    \left[
        \eta_j(1-e^{-v_j})-s_j v_j
    \right]
    =
    s_j\log\left(\frac{s_j}{\eta_j}\right)-s_j+\eta_j.
\]
Therefore, by convexity, 
\[
\begin{aligned}
    \frac{1}{4\delta^2}\WFR_\delta^2(\rho,\eta)
    =
    \inf_{\gamma\ge0}
    \left\{
        \sum_{i,j}c_{ij}\gamma_{ij}
        +
        \KL(\gamma\one\,|\,\rho)
        +
        \KL(\gamma^T\one\,|\,\eta)
    \right\}.
\end{aligned}
\]
This is exactly \eqref{eq:et-wfr}.
\end{proof}

\subsection{\textbf{Proof of Proposition~\ref{prop:et-kkt}}}
\label{app:proof-et-kkt}

\begin{proof}
Consider the entropy transport objective
\[
    \Phi(\gamma)
    =
    \sum_{i,j}c_{ij}\gamma_{ij}
    +
    \KL(\gamma\one\,|\,\rho)
    +
    \KL(\gamma^T\one\,|\,\eta).
\]
Set $I=\{i:\rho_i>0\}$ and $J=\{j:\eta_j>0\}$, both nonempty. Since $\Phi(0)=|\rho|+|\eta|<\infty$, the KL terms force $\gamma^*_{ij}=0$ outside $I\times J$. All indices and sums below are restricted to these supports. By extending the recovered $a^*,b^*$ by zero elsewhere, the omitted zero-mass marginal constraints automatically hold. Let
\[
    r_i=(\gamma\one)_i,\qquad s_j=(\gamma^T\one)_j.
\]
Since $\rho_i,\eta_j>0$ on $I\times J$ and $c_{ij}<\infty$, any minimizer satisfies $r_i^*>0$ for $i\in I$ and $s_j^*>0$ for $j\in J$. Because if \(r_i^\ast=0\), adding a small positive mass to some entry in the \(i\)-th row decreases the objective, since the derivative of \(r\log(r/\rho_i)\) at \(r=0\) is \(-\infty\). The same argument applies to \(s_j^\ast\).

For \(\gamma_{ij}^\ast>0\), the first-order optimality condition gives
\[
    0=
    c_{ij}
    +
    \log\left(\frac{r_i^\ast}{\rho_i}\right)
    +
    \log\left(\frac{s_j^\ast}{\eta_j}\right).
\]
For \(\gamma_{ij}^\ast=0\), the one-sided derivative must be nonnegative. Hence
\[
    c_{ij}
    +
    \log\left(\frac{r_i^\ast}{\rho_i}\right)
    +
    \log\left(\frac{s_j^\ast}{\eta_j}\right)
    \ge 0,
\]
and
\[
    \gamma_{ij}^\ast
    \left[
    c_{ij}
    +
    \log\left(\frac{r_i^\ast}{\rho_i}\right)
    +
    \log\left(\frac{s_j^\ast}{\eta_j}\right)
    \right]
    =0.
\]
These are the complementary slackness conditions on $I\times J$. Now we construct the original static WFR variables by
\[
    a_{ij}^\ast
    =
    \gamma_{ij}^\ast\frac{\rho_i}{r_i^\ast},
    \qquad
    b_{ij}^\ast
    =
    \gamma_{ij}^\ast\frac{\eta_j}{s_j^\ast}.
\]
They are feasible for the semi-coupling problem, since
\[
    \sum_j a_{ij}^\ast
    =
    \frac{\rho_i}{r_i^\ast}\sum_j\gamma_{ij}^\ast
    =
    \rho_i,
\]
and
\[
    \sum_i b_{ij}^\ast
    =
    \frac{\eta_j}{s_j^\ast}\sum_i\gamma_{ij}^\ast
    =
    \eta_j.
\]

It remains to show that this feasible pair is optimal for the original static WFR problem. If \(\gamma_{ij}^\ast>0\), then the complementary slackness condition gives
\[
    -2\log\omega_{ij}
    +
    \log\left(\frac{r_i^\ast}{\rho_i}\right)
    +
    \log\left(\frac{s_j^\ast}{\eta_j}\right)
    =
    0.
\]
Equivalently,
\[
    \frac{r_i^\ast s_j^\ast}{\rho_i\eta_j}
    =
    \omega_{ij}^2.
\]
Therefore,
\[
\begin{aligned}
    \omega_{ij}\sqrt{a_{ij}^\ast b_{ij}^\ast}
    &=
    \omega_{ij}\gamma_{ij}^\ast
    \sqrt{
        \frac{\rho_i\eta_j}{r_i^\ast s_j^\ast}
    }=
    \gamma_{ij}^\ast.
\end{aligned}
\]
If \(\gamma_{ij}^\ast=0\), then \(a_{ij}^\ast=b_{ij}^\ast=0\), so the same identity still holds:
\[
    \omega_{ij}\sqrt{a_{ij}^\ast b_{ij}^\ast}
    =
    \gamma_{ij}^\ast.
\]

Hence the static WFR objective at \((a^\ast,b^\ast)\) is
\[
\begin{aligned}
    &\sum_{i,j}
    \left[
        a_{ij}^\ast+b_{ij}^\ast
        -2\omega_{ij}\sqrt{a_{ij}^\ast b_{ij}^\ast}
    \right]=
    \sum_i\rho_i+\sum_j\eta_j
    -2\sum_{i,j}\gamma_{ij}^\ast.
\end{aligned}
\]
On the other hand, by the complementary slackness relation,
\[
\begin{aligned}
    \Phi(\gamma^\ast)
    &=
    \sum_{i,j}c_{ij}\gamma_{ij}^\ast
    +
    \sum_i
    \left[
        r_i^\ast\log\left(\frac{r_i^\ast}{\rho_i}\right)
        -r_i^\ast+\rho_i
    \right]+
    \sum_j
    \left[
        s_j^\ast\log\left(\frac{s_j^\ast}{\eta_j}\right)
        -s_j^\ast+\eta_j
    \right]\\
    &=
    \sum_i\rho_i+\sum_j\eta_j
    -2\sum_{i,j}\gamma_{ij}^\ast+
    \sum_{i,j}\gamma_{ij}^\ast
    \left[
        c_{ij}
        +
        \log\left(\frac{r_i^\ast}{\rho_i}\right)
        +
        \log\left(\frac{s_j^\ast}{\eta_j}\right)
    \right]\\
    &=
    \sum_i\rho_i+\sum_j\eta_j
    -2\sum_{i,j}\gamma_{ij}^\ast.
\end{aligned}
\]
Thus the entropy transport objective at \(\gamma^\ast\) equals to the static WFR objective at the constructed pair \((a^\ast,b^\ast)\).

By Proposition~\ref{prop:entropy}, the entropy transport problem and the static WFR semi-coupling problem have the same optimal value. Since \(\gamma^\ast\) is optimal for the entropy transport problem and \((a^\ast,b^\ast)\) is feasible for the static WFR problem with the same objective value, \((a^\ast,b^\ast)\) is an optimal solution of the original static WFR problem.
\end{proof}

\subsection{\textbf{Proof of Theorem~\ref{thm:udfm_cudfm_equivalence}}}
\label{app:proof-udfm-equivalence}

\begin{proof}
Let $\mathbb E_t$ denote $\mathbb E_{t\sim U[0,1]}$. Expanding the
marginal objective \eqref{eq:udfm_loss} gives
\begin{align*}
\mathcal L_{\mathrm{UDFM}}(\theta,\phi)
&=
\mathbb E_t
\sum_{x\in\mathcal X} p_t(x)
\Bigl[
(1-t)^2\big\|q_t^\theta(x,\cdot)-q_t(x,\cdot)\big\|_{2,x}^2
+
\lambda
\big|g_t^\phi(x)-g_t(x)\big|^2
\Bigr]
\\
&=
\mathbb E_t
\sum_{x\in\mathcal X} p_t(x)
\left[
(1-t)^2\Big(\big\|q_t^\theta(x,\cdot)\big\|_{2,x}^2
-
2\big\langle q_t^\theta(x,\cdot),q_t(x,\cdot)\big\rangle_x
+
\big\|q_t(x,\cdot)\big\|_{2,x}^2\Big)
\right.
\\
&\hspace{3.6cm}
\left.
+
\lambda\big(g_t^\phi(x)\big)^2
-
2\lambda g_t^\phi(x)g_t(x)
+
\lambda\big(g_t(x)\big)^2
\right]
\\
&=
\mathbb E_t
\sum_{x\in\mathcal X} p_t(x)
\Bigl[(1-t)^2\Big(
\big\|q_t^\theta(x,\cdot)\big\|_{2,x}^2
-
2\big\langle q_t^\theta(x,\cdot),q_t(x,\cdot)\big\rangle_x\Big)
\\
&\hspace{2.7cm}
+
\lambda\big(g_t^\phi(x)\big)^2
-
2\lambda g_t^\phi(x)g_t(x)
\Bigr]
+
C_{\mathrm{marg}},
\end{align*}
where
\[
C_{\mathrm{marg}}
=
\mathbb E_t
\sum_{x\in\mathcal X} p_t(x)
\left[
(1-t)^2\big\|q_t(x,\cdot)\big\|_{2,x}^2
+
\lambda \big(g_t(x)\big)^2
\right].
\]
is independent of $(\theta,\phi)$. By the marginalization identities
\eqref{eq:udfm_marginal_path}, \eqref{eq:udfm_marginal_rate}, and
\eqref{eq:udfm_marginal_growth}, we have
\[
p_t(x)
=
\int_{\mathcal Z}p_t(x\mid z)\,\pi(z)\dd z,
\]
\[
p_t(x)q_t(x,y)
=
\int_{\mathcal Z}
p_t(x\mid z)q_t(x,y\mid z)\,\pi(z)\dd z,
\qquad y\neq x,
\]
and
\[
p_t(x)g_t(x)
=
\int_{\mathcal Z}
p_t(x\mid z)g_t(x\mid z)\,\pi(z)\dd z.
\]
Substituting these identities into the parameter-dependent terms above gives
\begin{align*}
\mathcal L_{\mathrm{UDFM}}(\theta,\phi)
&=
\mathbb E_t
\int_{\mathcal Z}
\sum_{x\in\mathcal X}p_t(x\mid z)
\Bigl[
(1-t)^2\Big(\big\|q_t^\theta(x,\cdot)\big\|_{2,x}^2
-
2\big\langle q_t^\theta(x,\cdot),q_t(x,\cdot\mid z)\big\rangle_x\Big)
\\
&\hspace{2.2cm}
+
\lambda\big(g_t^\phi(x)\big)^2
-
2\lambda g_t^\phi(x)g_t(x\mid z)
\Bigr]\pi(z)\dd z
+
C_{\mathrm{marg}}.
\end{align*}
Next, we add and subtract the conditional target norm. Specifically,
\begin{align*}
&\big\|q_t^\theta(x,\cdot)\big\|_{2,x}^2
-
2\big\langle q_t^\theta(x,\cdot),q_t(x,\cdot\mid z)\big\rangle_x
\\
&\qquad =
\big\|
q_t^\theta(x,\cdot)-q_t(x,\cdot\mid z)
\big\|_{2,x}^2
-
\big\|q_t(x,\cdot\mid z)\big\|_{2,x}^2,
\end{align*}
and
\[
\lambda\big(g_t^\phi(x)\big)^2
-
2\lambda g_t^\phi(x)g_t(x\mid z)
=
\lambda\big|g_t^\phi(x)-g_t(x\mid z)\big|^2
-
\lambda\big(g_t(x\mid z)\big)^2.
\]
Therefore,
\begin{align*}
\mathcal L_{\mathrm{UDFM}}(\theta,\phi)
&=
\mathbb E_t
\int_{\mathcal Z}
\sum_{x\in\mathcal X}p_t(x\mid z)
\Bigl[
(1-t)^2\big\|
q_t^\theta(x,\cdot)-q_t(x,\cdot\mid z)
\big\|_{2,x}^2\\
&+
\lambda
\big|
g_t^\phi(x)-g_t(x\mid z)
\big|^2
\Bigr]\pi(z)\dd z
+C_{\mathrm{marg}}
\\
&-
\mathbb E_t
\int_{\mathcal Z}
\sum_{x\in\mathcal X}p_t(x\mid z)
\Bigl[
(1-t)^2\big\|q_t(x,\cdot\mid z)\big\|_{2,x}^2
+
\lambda\big(g_t(x\mid z)\big)^2
\Bigr]\pi(z)\dd z.
\end{align*}
Define
\[
C_{\mathrm{cond}}
=
\mathbb E_t
\int_{\mathcal Z}
\sum_{x\in\mathcal X}p_t(x\mid z)
\left[
(1-t)^2\big\|q_t(x,\cdot\mid z)\big\|_{2,x}^2
+
\lambda\big(g_t(x\mid z)\big)^2
\right]\pi(z)\dd z.
\]
The first term in the previous display is precisely
$\mathcal L_{\mathrm{CUDFM}}(\theta,\phi)$ in \eqref{eq:cudfm_loss}. Hence
\[
\mathcal L_{\mathrm{UDFM}}(\theta,\phi)
=
\mathcal L_{\mathrm{CUDFM}}(\theta,\phi)
+
C_{\mathrm{marg}}
-
C_{\mathrm{cond}}.
\]
Setting
\[
C:=C_{\mathrm{marg}}-C_{\mathrm{cond}},
\]
we obtain
\[
\mathcal L_{\mathrm{UDFM}}(\theta,\phi)
=
\mathcal L_{\mathrm{CUDFM}}(\theta,\phi)
+
C.
\]
Both $C_{\mathrm{marg}}$ and $C_{\mathrm{cond}}$ depend only on the target
marginal and conditional fields, and therefore $C$ is independent of
$(\theta,\phi)$. Differentiating both sides with respect to $(\theta,\phi)$
gives
\[
\nabla_{\theta,\phi}\mathcal L_{\mathrm{UDFM}}(\theta,\phi)
=
\nabla_{\theta,\phi}\mathcal L_{\mathrm{CUDFM}}(\theta,\phi).
\]
This proves the claim.
\end{proof}



\subsection{\textbf{Proof of Proposition~\ref{prop:approx_coarse}}}
\label{app:proof-total-mass-stability}

\begin{proof}
Let \(\zeta^0_{ij}=\bar\zeta_{P(i)P(j)}\) be the block-homogeneous lifting of the coarse affinity. For any feasible $(a,b)$,
\[
\begin{aligned}
|\mathcal J_\zeta(a,b)-\mathcal J_{\zeta^0}(a,b)|
\leq
\frac{\varepsilon}{4\delta^2}\sum_{i,j}\sqrt{a_{ij}b_{ij}}
\leq
\frac{\varepsilon}{4\delta^2}
\sqrt{\sum_{i,j}a_{ij}}
\sqrt{\sum_{i,j}b_{ij}}
=
\frac{\varepsilon}{4\delta^2}\sqrt{RS}.
\end{aligned}
\]
Hence the corresponding optimal values differ by at most the same amount. By Theorem~\ref{thm:exact_coarse}, the optimal value associated with $\zeta^0$ equals the coarse optimal value, which proves
\eqref{eq:wfr cost consistency} after multiplying by $4\delta^2$. 

For any semi-coupling, the travelling-Dirac formula \eqref{eq:travelling Dirac mass} gives the total mass of the induced measure flow
\[
M_f(t)
=
(1-t)^2R+t^2S
+
2t(1-t)
\sum_{i,j}\zeta_{ij}\sqrt{a_{ij}b_{ij}}.
\]
Since \(\mathcal J_\zeta(a,b)=R+S-2\sum_{i,j}\zeta_{ij}\sqrt{a_{ij}b_{ij}},\) an optimal fine WFR flow satisfies
\begin{equation}
M_f(t)
=
(1-t)R+tS
-
\frac{t(1-t)}{4\delta^2}
\mathrm{WFR}_{\delta}^2(\rho,\eta).    
\end{equation}
Analogously,
\begin{equation}
M_c(t)
=
(1-t)R+tS
-
\frac{t(1-t)}{4\delta^2}
\mathrm{WFR}_{\delta}^2(\bar\rho,\bar\eta).    
\end{equation}
Direct subtraction gives \eqref{eq:mass consistency}. Finally,
\begin{equation}
M_c(t)
\geq
(1-t)^2R+t^2S
\geq
2t(1-t)\sqrt{RS},
\end{equation}
which implies \eqref{eq:RME consistency}.
\end{proof}

\subsection{\textbf{Proof of Proposition~\ref{prop:path_stability}}}
\label{app:proof-measure-flow-convergence}

\begin{proof}
For each $n$, define
\begin{equation*}
\begin{aligned}
\zeta_{ij}^{(n)}
&=
\cos\left(
\frac{\sqrt{C_{ij}^{(n)}}}{2\delta}
\right),
\qquad
\zeta_{ij}^0
&=
\bar\zeta_{P(i)P(j)}
=
\cos\left(
\frac{\sqrt{\bar C_{P(i)P(j)}}}{2\delta}
\right).
\end{aligned}
\end{equation*}
By the same estimate used in
Proposition~\ref{prop:approx_coarse},
\begin{equation}
\label{eq:zeta_n_bound}
\max_{i,j}
\left|
\zeta_{ij}^{(n)}
-
\zeta_{ij}^0
\right|
\leq
\frac{\varepsilon_n}{8\delta^2}.
\end{equation}

Let $\Gamma(\rho,\eta)$ denote the fine-scale feasible set of semi-couplings,
\[
\Gamma(\rho,\eta)
=
\left\{
(a,b)\geq0:
\sum_j a_{ij}=\rho_i,
\sum_i b_{ij}=\eta_j
\right\}.
\]
Since
\[
0\leq a_{ij}\leq\rho_i,
\qquad
0\leq b_{ij}\leq\eta_j,
\]
the set $\Gamma(\rho,\eta)$ is closed and bounded in a finite-dimensional
Euclidean space, therefore compact. Consider any subsequence of
$\{(a^{(n)},b^{(n)})\}_{n\geq1}$.
By compactness, it contains a further subsequence, still indexed by
$n$ for simplicity, such that
\begin{equation}
\label{eq:semicoupling convergence}
(a^{(n)},b^{(n)})
\longrightarrow
(a^*,b^*).
\end{equation}
We first show that $(a^*,b^*)$ is optimal for the limiting
block-homogeneous fine-scale problem associated with
$\zeta^0$. For every $(a,b)\in\Gamma(\rho,\eta)$, \eqref{eq:zeta_n_bound} and
Cauchy--Schwarz imply
\begin{equation}
\label{eq:uniform_objective}
\begin{aligned}
\left|
\mathcal J_{\zeta^{(n)}}(a,b)
-
\mathcal J_{\zeta^0}(a,b)
\right|
&\leq
2
\max_{i,j}
\left|
\zeta_{ij}^{(n)}-\zeta_{ij}^0
\right|
\sum_{i,j}\sqrt{a_{ij}b_{ij}}\\
&\leq
\frac{\varepsilon_n}{4\delta^2}
\sqrt{
\sum_{i,j}a_{ij}
}
\sqrt{
\sum_{i,j}b_{ij}
}\\
&=
\frac{\varepsilon_n}{4\delta^2}\sqrt{RS}.
\end{aligned}
\end{equation}
Importantly, the right-hand side is independent of $(a,b)$.
Hence
$\mathcal J_{\zeta^{(n)}}$ converges uniformly to
$\mathcal J_{\zeta^0}$ on $\Gamma(\rho,\eta)$.

Now fix an arbitrary $(a,b)\in\Gamma(\rho,\eta)$.
Since $(a^{(n)},b^{(n)})$ is optimal for
$\mathcal J_{\zeta^{(n)}}$,
\[
\mathcal J_{\zeta^{(n)}}(a^{(n)},b^{(n)})
\leq
\mathcal J_{\zeta^{(n)}}(a,b).
\]
Using \eqref{eq:uniform_objective},
\[
\mathcal J_{\zeta^0}(a^{(n)},b^{(n)})
\leq
\mathcal J_{\zeta^0}(a,b)
+
\frac{\varepsilon_n}{2\delta^2}\sqrt{RS}.
\]
Taking $n\to\infty$ and using
\eqref{eq:semicoupling convergence} together with the continuity of
$\mathcal J_{\zeta^0}$ gives \(
\mathcal J_{\zeta^0}(a^*,b^*)
\leq
\mathcal J_{\zeta^0}(a,b).
\) Since $(a,b)\in\Gamma(\rho,\eta)$ was arbitrary,
$(a^*,b^*)$ is an optimal semi-coupling for the limiting block-homogeneous fine-scale problem.

Let
\[
\bar a^*_{AB}
=
\sum_{i\in B_A,j\in B_B}a^*_{ij},
\qquad
\bar b^*_{AB}
=
\sum_{i\in B_A,j\in B_B}b^*_{ij}.
\]
Since the limiting cost satisfies
\[
C^0_{ij}
=
\bar C_{P(i)P(j)},
\]
Theorem~\ref{thm:exact_coarse} implies that
$(\bar a^*,\bar b^*)$ is an optimal coarse-scale semi-coupling.
Moreover, if $\mu_t^*$ denotes the fine-scale WFR flow reconstructed
from $(a^*,b^*)$ under the limiting block-homogeneous cost, then
\begin{equation}
\label{eq:limiting_path_lift}
P_\#\mu_t^*
=
\bar\mu_t^*,
\end{equation}
where $\bar\mu_t^*$ is an optimal coarse-scale WFR flow generated by $(\bar a^*,\bar b^*)$.

It remains to show that the dynamic flows themselves
converge. For every pair $(i,j)$, let
\[
\theta_{ij}^{(n)}
=
\frac{\sqrt{C_{ij}^{(n)}}}{2\delta},
\qquad
\theta_{ij}^{0}
=
\frac{\sqrt{\bar C_{P(i)P(j)}}}{2\delta}.
\]
By \eqref{eq:epsilon_n},
\[
\theta_{ij}^{(n)}
\longrightarrow
\theta_{ij}^{0}.
\]
The travelling Dirac associated with
$(a_{ij}^{(n)},b_{ij}^{(n)},\theta_{ij}^{(n)})$ has mass
\[
m_{ij}^{(n)}(t)
=
(1-t)^2a_{ij}^{(n)}
+
t^2b_{ij}^{(n)}
+
2t(1-t)
\sqrt{a_{ij}^{(n)}b_{ij}^{(n)}}
\cos\theta_{ij}^{(n)}.
\]
Therefore,
\[
m_{ij}^{(n)}(t)
\longrightarrow
m_{ij}^*(t)
\]
uniformly for $t\in[0,1]$.

The corresponding progress variable can be represented through the
complex-plane formulation of the travelling Dirac:
\[
z_{ij}^{(n)}(t)
=
(1-t)\sqrt{a_{ij}^{(n)}}
+
t\sqrt{b_{ij}^{(n)}}
e^{\mathrm i\theta_{ij}^{(n)}},
\]
with
\[
m_{ij}^{(n)}(t)=|z_{ij}^{(n)}(t)|^2,
\qquad
k_{ij}^{(n)}(t)
=
\frac{\arg z_{ij}^{(n)}(t)}
{\theta_{ij}^{(n)}},
\]
with the usual continuous extension when the angle or an endpoint
mass is zero.
By
\eqref{eq:semicoupling convergence} and
$\theta_{ij}^{(n)}\to\theta_{ij}^0$,
the curves $z_{ij}^{(n)}$ converge uniformly to $z_{ij}^*$.

Away from points where $z_{ij}^*(t)=0$, continuity of the argument
implies
\[
k_{ij}^{(n)}(t)\to k_{ij}^*(t).
\]
At a point where $z_{ij}^*(t)=0$, although the progress variable
itself need not be uniquely defined, the two measure components
satisfy
\[
0\leq
m_{ij}^{(n)}(t)(1-k_{ij}^{(n)}(t))
\leq m_{ij}^{(n)}(t),
\qquad
0\leq
m_{ij}^{(n)}(t)k_{ij}^{(n)}(t)
\leq m_{ij}^{(n)}(t),
\]
and both converge to zero.
Hence the measure-valued travelling Dirac
\[
\mu_t^{ij,(n)}
=
m_{ij}^{(n)}(t)
\left[
(1-k_{ij}^{(n)}(t))\delta_i
+
k_{ij}^{(n)}(t)\delta_j
\right]
\]
depends continuously on its endpoint masses and cost, including at
zero endpoint masses. Consequently,
\[
\sup_{t\in[0,1]}
\left\|
\mu_t^{ij,(n)}
-
\mu_t^{ij,*}
\right\|_1
\longrightarrow0.
\]
Since the state space is finite, summing over all $(i,j)$ gives
\begin{equation}
\label{eq:fine_path_convergence}
\sup_{t\in[0,1]}
\left\|
\mu_t^{(n)}
-
\mu_t^*
\right\|_1
\longrightarrow0.
\end{equation}

Note that the pushforward $P_\#$ is linear and nonexpansive in the
$\ell^1$ norm.
\[
\|P_\#\mu-P_\#\nu\|_1
\leq
\|\mu-\nu\|_1.
\]
Combining this with
\eqref{eq:fine_path_convergence} and
\eqref{eq:limiting_path_lift} yields
\[
\sup_{t\in[0,1]}
\left\|
P_\#\mu_t^{(n)}
-
\bar\mu_t^*
\right\|_1
\longrightarrow0
\]
along the chosen subsequence. Therefore, every subsequential limit of the aggregated fine-scale
flows is an optimal coarse-scale WFR flow.

Now suppose that the coarse optimal flow is unique.
If \eqref{eq:path consistency} were false, there would exist
$c>0$ and a subsequence such that
\[
\sup_{t\in[0,1]}
\left\|
P_\#\mu_t^{(n)}
-
\bar\mu_t^*
\right\|_1
\geq c
\]
for every element of this subsequence. By compactness of the feasible semi-coupling set, this subsequence contains a further subsequence to which the preceding argument applies. Its aggregated flows must converge to an optimal coarse flow, which by uniqueness is $\bar\mu_t^*$. This contradicts the above lower bound. Hence \eqref{eq:path consistency} holds for the entire sequence.

Finally, assume $R,S>0$.
For every coarse optimal WFR flow,
\[
|\bar\mu_t^*|
\geq
(1-t)^2R+t^2S
\geq
\frac{RS}{R+S}>0,
\qquad t\in[0,1].
\]
For any two nonzero nonnegative measures $\mu,\nu$,
\[
\mathrm{TV}
\left(
\frac{\mu}{|\mu|},
\frac{\nu}{|\nu|}
\right)
=
\frac12
\left\|
\frac{\mu}{|\mu|}
-
\frac{\nu}{|\nu|}
\right\|_1
\leq
\frac{\|\mu-\nu\|_1}{|\nu|}.
\]
Taking
\[
\mu=P_\#\mu_t^{(n)},
\qquad
\nu=\bar\mu_t^*,
\]
we obtain
\[
\mathrm{TV}
\left(
\frac{P_\#\mu_t^{(n)}}{|P_\#\mu_t^{(n)}|},
\frac{\bar\mu_t^*}{|\bar\mu_t^*|}
\right)
\leq
\frac{R+S}{RS}
\left\|
P_\#\mu_t^{(n)}
-
\bar\mu_t^*
\right\|_1.
\]
Taking the supremum over $t$ and using
\eqref{eq:path consistency} proves
\eqref{eq:TV consistency}.
\end{proof}

\clearpage

\bibliographystyle{siamplain}
\bibliography{references}

\end{document}